\nonstopmode

\documentclass[]{amsart}
\usepackage{amssymb}
\usepackage{amscd}
\usepackage{graphicx}
\usepackage{epsfig}
\usepackage[all]{xy}
\newdir{ >}{{}*!/-10pt/@{>}}
\newdir^{ (}{{}*!/-5pt/@^{(}}
\newdir_{ (}{{}*!/-5pt/@_{(}}
\def\cgaps#1{}
\def\Cgaps#1{}
\allowdisplaybreaks

\def\undersetbrace#1\to#2{\underbrace{#2}_{#1}}
\def\oversetbrace#1\to#2{\overbrace{#2}^{#1}}
\def\AMSunderset#1\to#2{\underset{#1}{#2}}
\def\AMSoverset#1\to#2{\overset{#1}{#2}}

\swapnumbers

\newtheorem{prop}[subsection]{Proposition}
\newtheorem*{prop*}{Proposition}
\newtheorem{thm}[subsection]{Theorem}
\newtheorem*{thm*}{Theorem}

\newtheorem*{lem*}{Lemma}

\newtheorem*{cor*}{Corollary}
\newenvironment{demo}[1]{\par\smallskip\noindent{\bf #1.}}{\par\smallskip}
\newcommand{\scp}[2]{{\left\langle {#1}\, , \, {#2}\right\rangle}}
\newcommand{\avg}[1]{\overline{\left\langle {#1} \right\rangle}}

\def\cit#1#2{\ifx#1!\cite{#2}\else#2\fi} 
\def\ign#1{}             
\def\o{\circ}

\def\al{\alpha}
\def\be{\beta}
\def\ga{\gamma}
\def\de{\delta}

\def\et{\eta}
\def\th{\theta}

\def\ka{\kappa}
\def\la{\lambda}

\def\si{\sigma}

\def\ph{\varphi}

\def\ps{\psi}

\def\Ph{\Phi}
\def\Ps{\Psi}

\def\i{^{-1}}
\def\x{\times}
\def\p{\partial}
\let\on=\operatorname
\def\onb#1{\operatorname{\bf #1}}

\def\AMSonly#1{}

\def\op{\text{op}}
\def\od{{\text{od}}}
\def\ev{{\text{ev}}}
\begin{document}

\title
{A Metric on Shape Space with Explicit Geodesics}
\author{Laurent Younes, Peter W. Michor,
Jayant Shah, David Mumford}
\address{
Laurent Younes:
Johns Hopkins University,
Baltimore, MD 21218}
\email{laurent.younes@jhu.edu}

\address{
Peter W. Michor:
Fakult\"at f\"ur Mathematik, Universit\"at Wien,
Nordbergstrasse 15, A-1090 Wien, Austria; {\it and:}
Erwin Schr\"odinger Institut f\"ur Mathematische Physik,
Boltzmanngasse 9, A-1090 Wien, Austria
}
\email{Peter.Michor@univie.ac.at}

\address{
Jayant M. Shah:
Northeastern University
Department of Mathematics
360 Huntington Avenue,
Boston, MA 02115, USA
}
\email{shah@neu.edu}

\address{
David Mumford:
Division of Applied Mathematics, Brown University,
Box F, Providence, RI 02912, USA}
\email{David\_{}Mumford@brown.edu}

\thanks{All authors were supported by 
    NSF-Focused Research Group: 
    The geometry, mechanics, and statistics of the infinite dimensional
    shape manifolds.
    PWM was also supported by FWF Project P~17108.}
\subjclass{Primary 58B20, 58D15, 58E40}
\keywords{Shape space, diffeomorphism group, Riemannian metric, Stiefel manifold}
\date{\today}
\def\LaTeXonly{}

\begin{abstract}This paper studies a specific metric on plane curves that
has the property of being isometric to classical manifold (sphere,
complex projective, Stiefel, Grassmann) modulo change of
parametrization, each of these classical manifolds being associated to
specific qualifications of the space of curves (closed-open, modulo
rotation etc\ldots) Using these isometries, we are able to explicitely
describe the geodesics, first in the parametric case, then by modding
out the paremetrization and considering horizontal vectors. We also
compute the sectional curvature for these spaces, and show, in
particular, that the space of closed curves modulo rotation and change
of parameter has positive curvature. Experimental results that
explicitly compute minimizing geodesics between two closed curves are
finally provided
\end{abstract}

\maketitle

\section*{Introduction}

The definition and study of spaces of plane shapes has recently met a
large amount of interest \cite{bcgj95,gk91,ksm02,ky06,you98,sm04}, 
and has important applications, in object recognition, for the analysis
of shape databases or in medical imaging. The theoretical background
involves the construction of infinite
dimensional manifolds of shapes \cite{ksm02, sm04}. The Riemannian
framework, in particular, is appealing, because it provides shape
spaces with a rich structure which is also useful for applications. 
A general discussion of several
classes of metrics that can be introduced for this purpose can be
found in \cite{MM3}. 

The present paper focuses on a particular Riemannian metric that has
very specific properties. This metric, which will 
be described in the next section, can be seen as a limit case of one
of the classes studied in \cite{MM3}, and would receive the label
$H_{1, \infty}$ in the nomenclature introduced therein. One of its
surprising properties is that it can be characterized as a image of
a Grassmann manifold by a suitably chosen Riemannian submersion. 
A consequence of this is the possibility to derive explicit 
geodesics in this shape space.

A precursor of the $H_{1, \infty}$ metric has been introduced in
\cite{you98, you99} and studied in the context of open plane
curves. It has also recently been used in \cite{msj05}.  Because the
metric is placed on curves modulo changes of parametrization, the
computation of geodesics naturally provides an elastic matching
algorithm.

The paper is organized as follows. We first provide the definitions
and notation that we will use for spaces of curves, the $H_{1,
\infty}$ metric and the classical manifolds that will be shown to be
isometric to it. We then study some local properties of the obtained
manifold, discussing in particular its geodesics and sectional
curvature. We finally provide experimental results for the numerical
computation of geodesics and the solution of the related elastic
matching problem.

\section{Spaces of Curves}
Throughout this paper, we will assume our plane curves are curves in
the {\it complex} plane $\mathbb C$. Then real inner products and 
$2 \times 2$ determinants of real 2-vectors are given by 
$\langle x,y\rangle = \on{Re}(\bar x y)$ and 
$\det(x,y)=\on{Im}(\bar x y)$.

We first recall the notations for various spaces of plane curves 
which we will need, some of which were introduced in the previous 
paper \cite{MM3}. For all questions about infinite dimensional analysis and
differential geometry we refer to \cite{km93}.  By
$$\on{Imm}_{\op}=\on{Imm}([0, 2\pi],\mathbb C)$$ 
we denote the space of $C^\infty$-immersions 
$c:[0, 2\pi] \rightarrow \mathbb C$. 
Here `$\op$' stands for {\it open curve}. $B_{i,\op}$ is the quotient of 
$\on{Imm}_\op$ by the group $\text{Diff}^+([0,2\pi])$ of 
$C^\infty$ increasing diffeomorphisms of $[0,2\pi]$. Next
$$\on{Imm}_{\ev}(S^1,\mathbb C),\quad \on{Imm}_{\od}(S^1,\mathbb C)$$ 
are the spaces of $C^\infty$-immersions
$c:S^1 \rightarrow \mathbb C$ of {\it even}, respectively {\it odd}
rotation degree. Here, $S^1$ is the unit circle in $\mathbb C$, which
will be identified in this paper to $\mathbb R/(2\pi\mathbb Z)$.
Then $B_{i,\ev}$ and $B_{i,\od}$ are the
quotients of $\on{Imm}_\ev$, respectively $\on{Imm}_\od$ by 
the group $\text{Diff}^+(S^1)$ of $C^\infty$ orientation
preserving diffeomorphisms of $S^1$. For example, $B_{i,\od}$ 
contains the simple closed plane curves, since they have index +1
or $-1$ (depending on how they are oriented). These are the main
focus of this study. To save us from enumerating 
special cases, we will often consider open curves as defined on 
$S^1$ but with a possible discontinuity at 0. We will also consider the 
quotients of these spaces by the group of translations, by the group 
of translations and rotations and the group of translations, rotations 
and scalings.

Using the notation of \cite{MM3}, we can introduce the basic metric
studied in this paper on these three spaces of immersions, but {\it modulo
translations}, as follows. 
Identify $T_c (\on{Imm}/\text{transl})$ with the set of vector
fields $h:S^1 \rightarrow \mathbb C$ along $c$ modulo constant vector
fields. Then we consider the limiting case of the scale invariant 
metric of Sobolev order 1 from \cite{MM3},~4.8:
\begin{equation}
\label{eq:gc}
G_c(h,h)=G^{\mathrm{imm,scal},1,\infty}_c(h,h) 
= \frac{1}{\ell(c)}\int_{S^1} |D_sh|^2.ds
\end{equation}
where, as in \cite{MM3}, $ds=|c_\th|\,d\th$ is arclength measure, 
$D_s=D_{s,c}=|c_\th|\i \p_\th$ is the derivative with respect to 
arc length, and $\ell(c)$ is the length of $c$. We also recall for
later use the notation $v=c_\th/|c_\th|$ for the unit tangent 
vector, and, as multiplication by $i$ is rotation by 90 degrees,
$n=i.v$ for the unit normal. Note that this 
metric is invariant with respect to reparametrizations of the curve 
$c$, hence it induces a metric which we also call $G$ on the 
quotient spaces $B_{i,\op}, B_{i,\ev}$ and $B_{i,\od}$ also modulo 
translations.

The geodesic equation in all these metrics is a simple limiting 
case of those worked out in \cite{MM3}. Suppose $c(\th,t)$ is a
geodesic. Then:
$$c_{tt} = D_s^{-1}\left( 
\langle D_sc_t,v \rangle D_sc_t - \tfrac12 |D_sc_t|^2 v \rangle \right) -
\avg{\langle D_sc_t,v \rangle}.c_t -
\frac{1}{2} \avg{|D_sc_t|^2}.D_s^{-2}(\ka.n)$$
Here the bar indicates the {\it average} of the quantity over the 
curve $c$, i.e.\ $\avg{F}=\tfrac{1}{\ell}\int F ds$. Unfortunately, this case was not worked out in \cite{MM3}, 
hence we give the details of its derivation in Appendix I.
The local existence and uniqueness of solutions to this equations
can be proved easily, essentially because of the regularizing
influence of the term $D_s^{-1}$. This will also follow from the explicit
representation of these geodesics to be given below, but because of
its more general applicability, we give a direct proof in Appendix I.

It is convenient to introduce the {\it momentum} $u=-D_s^2(c_t)$ 
associated to a geodesic. Using the momentum, the geodesic equation 
is easily rewritten in the more compact form:
$$ u_t = -\langle u, D_sc_t \rangle v - \left( \langle D_sc_t,v \rangle
-\avg{\langle D_sc_t,v \rangle}\right) u -\frac{1}{2} \left(
|D_sc_t|^2 + \avg{|D_sc_t|^2} \right) \ka(c).n.$$

By the theory of Riemannian submersions, geodesics on the quotient 
spaces $B_i$ are nothing more than {\it horizontal} geodesics in $\on{Imm}$, 
that is geodesics which are perpendicular at one hence all points to the
orbit of the group of reparametrizations. As is shown in \cite{MM3},
horizontality is equivalent to the condition $u=a.n$ for some scalar
function $a(\th,t)$. Substituting $u=a.n$ and taking the $n$-component
of the last equation, we find that horizontal geodesics are given by:
$$ a_t  = -a\left( 
\langle D_sc_t,v \rangle - \avg{\langle D_sc_t,v \rangle} \right) +
\frac{\ka(c)}{2} \left( |D_sc_t|^2 + \avg{|D_sc_t|^2} \right). $$

There are several conserved momenta along each geodesic 
$t\mapsto c(\th,t)$ of this metric (see \cite{MM3},~4.8):
The `reparametrization' momentum is
$$\frac{-1}{\ell(c)}\langle c_\th, D_{s,c}^2c_t \rangle |c_\th|$$
which vanishes along all horizontal geodesics. 
The translation momentum vanishes because the metric does not 
feel translation (constant vector fields along $c$). 
The angular momentum is 
$$\frac{-1}{\ell(c)}\int_{S^1}\langle i.c, D^2_{s,c}c_t \rangle\,ds
= \frac{1}{\ell(c)}\int_{S^1}\ka \langle v, c_t \rangle\,ds.$$
Since the metric invariant under scalings, we also have the scaling
momentum
$$\frac{-1}{\ell(c)}\int_{S^1} \langle c,D_{s,c}^2c_t \rangle\,ds 
= \p_t \log\ell(t).$$
So we may equivalently consider either the quotient space 
$\on{Imm}/\text{translations}$ or consider the section of the 
translation action $\{c\in \on{Imm}:c(0)=0\}$. In the same way, 
we may either pass modulo scalings or consider the section by 
fixing $\ell(c)=1$, since the scaling momentum vanishes here. 
Finally, in some cases, we will pass modulo rotations. We could 
consider the section where angular momentum vanishes: but this 
latter is not especially simple

\section{The Basic Mapping for Parametrized Curves}

\subsection{The basic mapping}
We introduce the three function spaces:
\begin{align*}
V_\op &= \text{Vector space of all } C^\infty \text{ mappings } 
f: [0,2\pi] \rightarrow \mathbb R \\
V_{\ev} &= \text{Vector space of all } C^\infty \text{ mappings } 
f: S^1 \rightarrow \mathbb R \text{ such that } f(\th+2\pi) 
\equiv f(\th)\\
V_{\od} &= \text{Vector space of all } C^\infty \text{ mappings } 
f: S^1 \rightarrow \mathbb R \text{ such that } f(\th+2\pi) 
\equiv -f(\th)
\end{align*}
All three spaces have the weak inner product: 
$$ \|f\|^2 = \int_0^{2\pi} f(x)^2 dx.$$
Given $e,f$ from any of these spaces, the basic map is:
$$ \Ph: (e,f) \longmapsto c(\th)=(1/2) \int^\th_0 
\left(e(x)+if(x)\right)^2 dx.$$
The map $c$ so defined carries $[0,2\pi]$ or $S^1$ to $\mathbb C$. 
It need not be an immersion, however, because $e$ and $f$ might
vanish simultaneously. Define 
$$Z(e,f)=\{\th : e(\th)=f(\th)=0\}.$$
Then we get three maps:
$$ \{ (e,f)\in V \times V : Z(e,f)=\emptyset\} \longrightarrow 
\on{Imm}_x, \text{ for } V=V_\op, V_\ev, V_\od.$$ 

Looking separately at the three cases, define first the sphere 
$\onb{S}(V^2_\op)$ to be the set of $(e,f)\in V_{\op}^2$
such that $\|e\|^2 + \|f\|^2 = 2$. $\onb{S}^0(V^2_\op)$ is defined
as the subset where $Z(e,f)=\emptyset$. Then the magic of the
map $\Ph$ is shown by the following key fact \cite{you98}:

\begin{thm} $\Ph$ defines a map:
$$ \Ph: \onb{S}^0(V^2_\op) \longrightarrow \left\{ c \in
\on{Imm}_{\op}:  \ell(c)=1, c(0)=0  \right\} \cong \on{Imm}_{\op} \Big/ 
(\text{transl,scalings})$$
which is an isometric 2-fold covering, using the natural metric 
on $\onb{S}$ and the metric $G^{\mathrm{imm,scal},1,\infty}$ on 
$\on{Imm}_{\op}$.
\end{thm}

\begin{demo}{Proof} The mapping $\Ph$ is surjective: Given
$c\in \on{Imm}_{\text{\op}}$ with $c(0)=0$ and $\ell(c)=1$,
we write $c'(u)=r(u)e^{i\al(u)}$. Then we may choose 
$e(u)=\sqrt{2r(u)}\cos\frac{\al(u)}2$ and
$f(u)=\sqrt{2r(u)}\sin\frac{\al(u)}2$. Since 
$1=\ell(c)=\int_0^{2\pi}|c'(u)|\,du = \int_0^{2\pi}r(u)\,du$
we see that $\|e\|^2+\|f\|^2=\int_0^{2\pi}
2r(u)(\cos^2(\frac{\al(u)}2)+\sin^2(\frac{\al(u)}2))\,du = 2$.
The only choice here is the sign of the square root, i.e.\ 
$\Ph(-e,-f)=\Ph(e,f)$, thus $\Ph$ is 2:1.

To see that $\Ph$ is an isometry, let $\Ph(e,f)=c=x+iy$ and 
$\de c = \de x+i\de y$. Then the differential $D\Ph(e,f)$ is 
given by
\begin{equation}
\label{eq:d.phi}
D\Phi(e,f): (\de e, \de f) \mapsto \de c(\th) =  \int^\th (\de e + i\de f)(e+if) d\th.
\end{equation}
We have $ds = (1/2) |e+if|^2 d\th$. This implies first that $\ell(c) = (\|e\|^2 +
\|f\|^2)/2 = 1$ as required. Then:  
\begin{align*}
D_s(\de c) &= \frac{2(e+if)(\de e+i\de f)}{|e+if|^2}
\\
G_c(\de c, \de c) &= 
\frac{1}{2} \int_0^{2\pi} \left|D_s(\de c)\right|^2 ds
= \int_0^{2\pi} |\de e + i\de f|^2 d\th = \|(\de e, \de f)\|^2.
\end{align*}
\end{demo}

The dictionary between pairs $(e,f)$ and immersions $c$ connects many properties of 
each with those of the other. Curvature $\ka$ works out especially nicely. 
We list here some of the connections:
\begin{align*}
\frac{ds}{d\th} = |c_\th| &= \tfrac12(e^2+f^2)\\
v = D_s(c) &= \frac{(e+if)^2}{e^2+f^2}
\end{align*}
and if $W_\th(e,f)=ef_\th-fe_\th$ is the Wronskian, then:
\begin{align*}
v_\th &= \left(\frac{(e+if)^2}{e^2+f^2}\right)_\th = 2\frac{W_\th(e,f)}{e^2+f^2}iv, 
\quad \text{hence}\\
\ka &= 2\frac{W_\th(e,f)}{(e^2+f^2)^2}\quad\text{for the curvature of  }c.
\end{align*}

\subsection{Geodesics leaving the space of immersions}
Since geodesics on a sphere are always given by great circles,
this theorem gives us the first case of explicit geodesics
on spaces of curves in the metric of this paper. However, note
that great circles in the open part $\onb{S}^{0}$ are susceptible
to crossing the `bad' part $\onb{S}-\onb{S}^{0}$ somewhere. This
occurs if and only if there exists $\th$ such that $(e+if)(\th)$ and
$(\de e+ i\de f(\th))$ have identical complex arguments modulo $\pi$.
So we
find that our metric on $\on{Imm}$ is {\it incomplete}.

\begin{figure}[ht]
 \begin{center}
\includegraphics[width=0.7\textwidth]{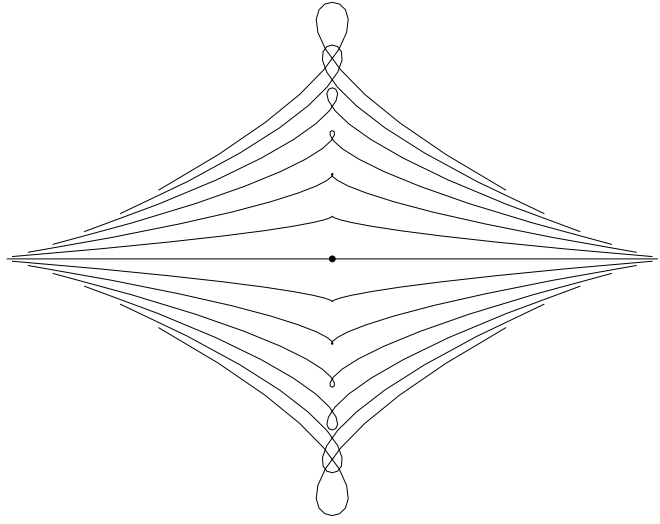}\\ 
\caption{\label{fig:genbifurc}\footnotesize The generic way in which a family 
of open immersions crosses the hypersurface where $Z \ne \emptyset$.
The parametrized straight line in the middle of the family has 
velocity with a double zero at the black dot, hence is not an
immersion. See text.}
\end{center}
\end{figure}

We can form a commutative diagram:
$$ \xymatrix{\Ph: \onb{S}^0(V^2_\op)\; 
\ar[r]^-{\text{2-fold}} \ar@{^{(}->}[d] &
\on{Imm}_{\op}\Big/ \left( \text{ transl,scalings }\right) \ar@{^{(}->}[d] \qquad \\  
\widetilde{\Ph}: \onb{S}(V_\op) \; 
\ar[r] &
C^{\infty}([0,2\pi],\mathbb C) \Big/ \left( \text{transl,scalings}\right)}
$$ 
where we have denoted the extended $\Ph$ by $\widetilde \Ph$. For 
rather technical reasons $\widetilde \Ph$ is not surjective: there
are pathological non-negative $C^\infty$ functions which have no 
$C^\infty$ square root, see \cite{klm04}, e.g. 
But what this diagram does do is give some 
space of maps to hold the extended geodesics. The example:
\begin{align*} 
e(x,s) +if(x,s) &= (x+is)/\sqrt{C},\quad \text{where } \\
&  -\pi \le x \le \pi, -1 \le s \le 1 
\text{ and } C=2a^3/3+2as^2\\
\text{and } c(x,s) &= (x^3/3-s^2x+isx^2)/C - is/2 \text{ (suitably translated)}
\end{align*}
is shown in figure \ref{fig:genbifurc}. This is a geodesic in which 
all curves are immersions for $s\ne 0$, but $c_x(x,0)$ has a double zero
at $x=\pi$.

\subsection{The basic mapping in the periodic case}
Next, consider the periodic cases. Here we need the {\it Stiefel
manifolds}: 
$$\onb{St}(2,V) = \{ \text{orthonormal pairs } (e,f) \in V\times V\}, V =
V_\ev \text{ or } V_\od$$
and $\onb{St}^0(2,V)$ the subset defined by the constraint  
$Z(e,f)=\emptyset$. For later use, it is also convenient to note that
$\onb{St}(2,V)=\{A\in L(\mathbb C,V): A^T.A= \on{Id}_{\mathbb C}\}$ 
where $L(\mathbb C,V)$ is the space of linear maps from $\mathbb C$ to
$V$ and $\langle A^\top v,w \rangle_{\mathbb C}=\langle v,Aw \rangle_V$. 
For $(e,f)\in \onb{S}^0(V^2_\op)$, when is $c=\Ph(e,f)$ periodic?
If and only if we have:
\begin{align*}
&c'=\tfrac12 (e+if)^2\quad\text{is periodic, so that }e,f\in
V_\ev \text{ or } V_\od;\tag{A} \\
&0 = \int_0^{2\pi}c'(u)\,du = \frac12 \int_0^{2\pi}(e^2-f^2)\,du  
+ i \int_0^{2\pi} ef\,du. \tag{B}
\end{align*}
Condition \thetag{B} says that $\|e\|^2=\|f\|^2=1$ (since the sum 
is 2) and $\langle e,f \rangle=0$, so that 
$(e,f)\in \onb{St}^0(2,V_\ev)$ or $(e,f)\in \onb{St}^0(2,V_\od)$.
Recall that the index $n$ of an immersed curve $c$ is defined by
considering $\log(c')$. The log must satisfy 
$\log(c'(\th+2\pi))\equiv \log(c'(\th))+2\pi n$ for some $n$ and 
this is the index. So this index is even or odd
depending on whether the square root of $c'$ is periodic or 
anti-periodic, that whether $(e,f)\in V_\ev$ or $\in V_\od$.
So if $\Ph$ is restricted to $\onb{St}^0(2, V_\ev)$ or 
$\onb{St}^0(2, V_\od)$ (and is still denoted $\Ph$), it provides 
isometric 2-fold coverings 
\begin{align*} 
\Ph: \onb{St}^0(2,V_\ev) &\longrightarrow 
\{c\in\on{Imm}_{\ev}(S^1, \mathbb C):c(0)=0, \ell(c)=1\}\\
\Ph: \onb{St}^0(2,V_\od) &\longrightarrow 
\{c\in\on{Imm}_{\od}(S^1, \mathbb C):c(0)=0, \ell(c)=1\}
\end{align*} 

All three of these maps $\Ph$ can be modified so as to divide out
by rotations. The mapping $(e,f) \mapsto e^{i\ph}(e+if)$ produces
a rotation of the immersed curve $\Ph(e,f)$ through an angle $2\ph$.
The complex projective space $\onb{CP}(V^2_\op)$ is $\onb{S}(V^2_{\op})$
divided by the action of rotations, and we denote by
$\onb{CP}^{0}(V^2_\op)$ the subset gotten by dividing 
$\onb{S}^{0}(V^2_{\op})$ by rotations. The group generated by
translations, rotations and scalings will be called the group
of {\it similitudes}, abreviated as `sim'. Then we get the variant:
$$ \Ph: \onb{CP}^{0}(V^2_{\op}) \longrightarrow \on{Imm}_\op /(\text{sim}).$$

Similarly, let $\onb{Gr}(2,V)$ be the Grassmannian of unoriented
2-dimensional subspaces of $V$ and let $\onb{Gr}^0(2,V)$ be 
the subset of those $W$ with $Z(W)=\emptyset$ for $V=V_\ev$ or $=V_\od$. 
Then we have maps:
\begin{align*} 
\Ph: \onb{Gr}^0(2,V_\ev) &\longrightarrow \on{Imm}_{\ev}/
(\text{sim})\\
\Ph: \onb{Gr}^0(2,V_\od) &\longrightarrow \on{Imm}_{\od}/
(\text{sim})
\end{align*} 

For later use, we describe the tangent spaces of these spaces.
The tangent space $T_W \onb{Gr}$ to $\onb{Gr}(2,V)$ at 
$W \in \onb{Gr}(2,V)$ is naturally identified with
$L(W,W^\bot)$ and has the following norm, induced from that
on $V$: 
$$ \|v\|^2 = \text{tr}(v^T \circ v) = \|v(e)\|^2 + \|v(f)\|^2$$ 
for $v \in T_W\onb{Gr}$ and $\{e,f\}$ an orthonormal basis of $W$. 
Similarly, $T_{(e,f)}\onb{St}$ can be naturally identified with pairs 
$\{\de e, \de f\}$ in $V$ such that $\langle e, \de e\rangle=
\langle f,\de f \rangle = \langle e,\de f\rangle + 
\langle f,\de e\rangle =0$ with norm: 
$$ \|(\de e, \de f)\|^2 = \|\de e\|^2 + \|\de f \|^2.$$
The same definition holds for  $T_{(e,f)}(\onb{S})$, this time with the
constraint  $\langle e, \de e\rangle+\langle f,\de f
\rangle =0$.

\section{The Basic Mapping for Shapes}

\subsection{Dividing out by the group of reparametrizations}
Let $C^{\infty,+}([0,2\pi])$ be the group of increasing
diffeomorphisms $\ph$ of $[0,2\pi]$ (so $\ph(0)=0, \ph(2\pi)=2\pi$) 
and let $C^{\infty,+}(\mathbb R)$ be the group of increasing
$C^\infty$ diffeomorphisms $\ph: \mathbb R \rightarrow \mathbb R$ 
such that $\ph(x+2\pi) \equiv \ph(x)+2\pi$ for all $x$. Modulo
the central subgroup of translations $\ph(x)=x+2\pi n$, the
second group is just $\onb{Diff}^+(S^1)$. For $V=V_\op, V_\ev, V_\od$
let $U(V)$ be the group of unitary maps on $V$ given by: 
$$ f \mapsto f^\ph = \sqrt{\ph'}.(f \circ \ph), 
\text{ where } 
\left\{ \begin{array}{l}\ph \in C^{\infty,+}([0,2\pi]) \text{ if }
V=V_{\op} \\
\ph \in C^{\infty,+}(\mathbb R) \text{ if } V=V_{\ev},V_{\od} \end{array} 
\right.$$
These are the reparametrization groups for our various spaces.
The infinitesimal action of a vector field $X$ on $[0,2\pi]$ or 
a periodic vector field $X$ on $S^1$ is then 
\begin{equation}
\label{eq:inf.act}
f\mapsto \tfrac12 X_\th.f + X.f_\th.
\end{equation}

For all three sets of isometries $\Ph$, we can now divide each side 
by the re\-pa\-ra\-metri\-za\-tion group $U(V)$. 
For open curves, we get a diagram
$$ \xymatrix{\overline{\Ph}: \onb{S}^0(V^2_\op)/U(V_\op)\; 
\ar[r]^-{\text{2-fold}} \ar[d] &
B_{i,\op}\Big/ \left( \text{ transl,scalings }\right) \ar[d] \quad\\  
\overline{\Ph}: \onb{CP}^0(V^2_\op)/U(V_\op) \; 
\ar[r]^-{\approx} &
B_{i,\op} \Big/ \left(\text{sim}\right)}
$$ 
and a similar one for closed curves of even and odd index where 
$V=V_\ev,V_\od$ and $B=B_{i,\ev}, B_{i,od}$:
$$ \xymatrix{\overline{\Ph}: \onb{St}^0(2,V)/U(V)\; 
\ar[r]^-{\text{2-fold}} \ar[d] &
B \Big/ \left( \text{transl,scalings}\right) \ar[d]\qquad\\  
\overline{\Ph}: \onb{Gr}^0(2,V)/U(V) \; 
\ar[r]^-\approx  &
B \Big/ \left(\text{sim}\right)}
$$ 

Here we have divided by isometries on both the left and right: 
by $U(V)$ or $U(V)\times S^1$ on the left (where $S^1$ rotates 
the basis $\{e,f\}$) and by reparametrizations and rotations 
on the right. Thus $\overline\Ph$ is again an isometry if we 
make both quotients into Riemannian submersions. This means 
we must identify the tangent spaces to the quotients with the 
{\it horizontal} subspaces of the tangent spaces in the larger 
space, i.e. those perpendicular to the orbits of the isometric 
group actions. For $\onb{St}$, this means: 

\begin{prop} \label{horizontalinSt}
The tangent vector $\{\de e, \de f\}$ to $\onb{St}$
satisfies:
$$ 
\langle \de e, e \rangle =\langle \de f,f \rangle = 0, \langle \de
e,f\rangle +\langle \de f, e \rangle =0.
$$
It is horizontal for the rotation action if and only if:
$$\text{both }\de e, \de f \text{ are perpendicular to both } e,f. $$
It is horizontal for the reparametrization group if:
$$ W_\th(e, \de e) + W_\th(f, \de f) = 0 $$
where $W_\th(a,b) = a.b_\th-b.a_\th$ is the Wronskian with respect to
the parameter $\th$.
\end{prop}
\begin{demo}{Proof}
Consider the action of rotations, which is one-dimensional, with
orbits $\be\mapsto e^{i\be}(e+if)$; the direction at $(e,f)$ is chosen as
$(-f, e)$. So $(\de e, \de f)$ being horizontal at $(e,f)$ for this
action means that $- \langle f, \de e\rangle + \langle e, \de f\rangle
= 0$. This proves the first assertion.

For the action of $U(V)$, one has to note that horizontal vectors must
satisfy 
$$
\langle \tfrac12 X_\th.e + X.e_\th , \de e\rangle + \langle \tfrac12 X_\th.f +
X.f_\th, \de f\rangle = 0
$$
for any periodic vector field $X$ on $\mathbb R$, which yields the
horizontality condition after integration by parts of the terms in
$X_\th$.
\qed\end{demo}

Horizontality on the shape space side means (see \cite{MM3}):

\begin{prop}
$h\in T_c\on{Imm}(S^1,\mathbb C)$ is horizontal for the action of
$\on{Diff}(S^1)$ if and only if
$D_s^2(h)$ is normal to the curve, i.e. $\langle v,D_s^2(h) \rangle =
0$.
\end{prop}

\begin{prop}
For any smooth path $c$ in $\on{Imm}(S^1,\mathbb R^2)$ there exists a
smooth
path $\ph$ in $\on{Diff}(S^1)$ with $\ph(0,\;.\;)=\on{Id}_{S^1}$ depending
smoothly on $c$ such that
the path $e$ given by $e(t,\th)=c(t,\ph(t,\th))$ is horizontal:
$\langle D_s^2(e_t),e_\th\rangle=0$.
\end{prop}
This is a variant of \cite[4.6]{MM3}.

\begin{demo}{Proof}
Writing $D_c$ instead of $D_s$ we note that 
$D_{c\o\ph}(f\o \ph)=\frac{(f_\th\o\ph)\ph_\th}{|c_\th\o\ph|.|\ph_\th|}
=(D_c(f))\o\ph$ for $\ph\in\on{Diff}^+(S^1)$. So we have 
$L_{n,c\o\ph}(f\o\ph)=(L_{n,c}f)\o\ph$.

Let us write $e=c\o \ph$ for $e(t,\th)=c(t,\ph(t,\th))$, etc.
We look for $\ph$
as the integral curve of a time dependent vector field $\xi(t,\th)$ on
$S^1$, given by $\ph_t=\xi\o \ph$.
We want the following expression to vanish:
\begin{align*}
\langle D_{c\o\ph}^2(\p_t(c\o\ph)),\p_\th(c\o\ph) \rangle
&=\langle D_{c\o\ph}^2(c_t\o\ph + (c_\th\o\ph)\,\ph_t),(c_\th\o\ph)\,\ph_\th
\rangle
\\&
=\langle D_{c}^2(c_t)\o\ph + D_{c}^2(c_\th.\xi)\o\ph,c_\th\o\ph
\rangle\ph_\th
\\&
=\bigl((\langle D_{c}^2(c_t),c_\th\rangle +\langle D_{c}^2(\xi.c_\th),
c_\th\rangle)\o\ph\bigr)\,\ph_\th.
\end{align*}
Using the time dependent vector field
$\xi=-\frac1{|c_\th|}D_c^{-2}(\langle D_{c}^{2}(c_t),v\rangle)$
and its flow $\ph$ achieves this.
\qed\end{demo}

\subsection{Bigger spaces}
\label{sec:bigger}
As we will see below, we can describe geodesics in the 
`classical' spaces $\onb{S}, \onb{CP}, \onb{St}, \onb{Gr}$ 
quite explicitly. By the above isometries, this gives 
us the geodesics in the various spaces $\on{Imm}, B_i$. 
BUT, as we mentioned above for the space $\onb{S}$, 
geodesics in the `good' parts $\onb{S}^0, \onb{CP}^0, 
\onb{St}^0, \onb{Gr}^0$ do not stay there, but they cross 
the `bad' part where $Z(e,f) \neq \emptyset$. Now 
the basic mapping is still defined on the full sphere,
projective space, Stiefel manifold or Grassmannian, 
giving us some smooth mappings of $[0,2\pi]$ or $S^1$ 
to $\mathbb C$, possibly modulo translations, rotations
and/or scalings. 

But when we divide by $U(V_\op)$, a major problem arises. 
The orbits of $U(V_\op)$ acting on $C^\infty([0,2\pi],\mathbb C)$
are not closed, hence the topological quotient of the space 
$C^\infty([0,2\pi],\mathbb C)$ by $U(V_\op)$ is {\it not 
Hausdorff}. This is shown by the following construction:
\begin{enumerate}
\item Start with a $C^\infty$ non-decreasing map $\ps$ from 
$[0,2\pi]$ to itself such that $\ps(\th)\equiv \pi$ for $\th$ 
in some interval $I$.
\item Let $\ps_n(\th)=(1-1/n).\ps(\th)+\th/n$. The sequence $\{\ps_n\}$ 
of diffeomorphisms of $[0,2\pi]$ converges to $\ps$.
\item Then for any $c \in \on{Imm}_\op$, the maps $c \circ \ps_n$ are
all in the orbit of $c$. But they converge to $c\circ \ps$ which 
is constant on the whole interval $I$, hence is not in the
orbit.
\end{enumerate}
Thus, if we want some Hausdorff space of curves which a) have 
singularities more complex than those of immersed curves and b) 
can hold the extensions of geodesics in some space $B_i$ which 
come from the map $\Ph$, we must divide $C^\infty([0,2\pi],
\mathbb C)$ by some equivalence relation larger than the group 
action by $U(V)$. The simplest seems to be: first define a
{\it monotone relation} $R \subset [0,2\pi]\times [0,2\pi]$ to be any
closed subset such that $p_1(R) = p_2(R) = [0, 2\pi]$ ($p_1$ and $p_2$
being the projections on the axes) and for every pair of points $(s_1,t_1) \in R$ and 
$(s_2,t_2) \in R$, either $s_1 \le s_2$ and $t_1 \le t_2$ or
vice versa. Then $f,g:[0,2\pi] \rightarrow \mathbb C$ are 
{\it Fr\'echet} equivalent if there is a monotone relation $R$
such that $f(s)=g(t), \forall (s,t) \in R$.

This is a good equivalence relation because if 
$\{f_n\},\{g_n\}:[0,2\pi] \rightarrow \mathbb C$ are two sequences 
and $\lim_n f_n = f, \lim_n g_n = g$ and $f_n,g_n$ are Fr\'echet
equivalent for all $n$, then $f,g$ are Fr\'echet equivalent. The
essential point is that the set of non-empty closed subsets of 
a compact metric space $X$ is compact in the Hausdorff topology
(see \cite{alexhopf}). Thus if $\{R_n\}$ are the monotone relations
instantiating the equivalence of $f_n$ and $g_n$, a subsequence
$\{R_{n_k}\}$ Hausdorff converges to some $R \subset [0,2\pi]
\times [0,2\pi]$ and it is immediate that $R$ is a monotone 
relation making $f$ and $g$ Fr\'echet equivalent.

Define 
$$B_{\text{big},\op} = C^\infty([0,2\pi],\mathbb C)/
\text{Fr\'echet equivalence, translations, scalings}.$$
Then we have a commutative diagram:
$$ \xymatrix{\Ph: \onb{S}^0(V^2_\op)\; 
\ar[r] \ar@{^{(}->}[d] &
B_{i,\op}\Big/ \left( \text{transl,scalings }\right) 
\ar@{^{(}->}[d] \qquad \\  
\widetilde{\Ph}: \onb{S}(V^2_\op)  \; \ar[r] & B_{big,\op}} $$ 
Thus the whole of a geodesic which enters the `bad' part of 
$\onb{S}(V_\op)$ creates a path in $B_{big,\op}$. Of course, the
same construction works for closed curves also. We will see
several examples in the next section.

\section{Construction of Geodesics}
\subsection{Great circles in Spheres} \label{greatcircles}
The space $\onb{S}(V^2)$ being the sphere of radius 
$\sqrt{2}$ on $V^2$, its geodesics are the great circles. 
Thus, the geodesic distance between $(e^0, f^0)$ and $(e^1, f^1)$ 
is given by $\sqrt{2}D$ with
$$
D  = \arccos\big( (\scp{e^0}{e^1} + \scp{f^0}{f^1})/2\big)
$$
and the geodesic is given by
\begin{eqnarray*}
e(t) &=& \frac{\sin((1-t)D)}{\sin D} e^0 + \frac{\sin (tD)}{\sin D} e^1\\
f(t) &=& \frac{\sin((1-t)D)}{\sin D} f^0 + \frac{\sin (tD)}{\sin D} f^1
\end{eqnarray*}
The corresponding geodesic on $\on{Imm}_{\op}$ modulo 
translation and scaling is the time-indexed family of curves 
$t\mapsto c(u, t)$ with
$$
\partial c/\partial u = \tfrac12(e(t)+if(t))^2 = (e(t)^2-f(t)^2)/2 + ie(t)f(t)
$$

The following notation will be used throughout this section: 
\begin{align*}
c^0(u) &=c(u,0),& c^1(u)&=c(u,1) \\  
\partial c^0 /\partial u &= r_0(u) e^{i\al^0(u)}, 
&\partial c^1 /\partial u &= r_1(u) e^{i\al^1(u)}
\end{align*} 
so that $e^j=\sqrt{2r_j}\cos\tfrac{\al^j}2$ and
$f^j=\sqrt{2r_j}\sin\tfrac{\al^j}2$ for $j=0,1$.
Thus the distance $D$ is
$$
D_{\op}(c^0, c^1) = \arccos \int_0^{2\pi}  \sqrt{r_0r_1} \cos
\frac{\al^1-\al^0}{2} du.
$$
The metric on $\on{Imm}_{\op}$ modulo rotations is 
\begin{align*}
&D_{\text{op, rot}}(c^0, c^1) 
= \inf_\al \arccos \int_0^{2\pi}  \sqrt{r_0r_1} 
  \cos \frac{\al^1-\al^0-\al}{2} du
\\& 
=\arccos \sup_\al \int_0^{2\pi}  \sqrt{r_0r_1} 
\left(\cos \frac{\al^1-\al^0}{2}\cos\frac\al 2 +
\sin\frac{\al^1-\al^0}2\sin\frac{\al}2\right) du
\\&
=  \arccos \left( \Big(\int_0^{2\pi}  \sqrt{r_0r_1} \cos
\frac{\al^1-\al^0}{2} du\Big)^2 
 + \Big(\int_0^{2\pi}  \sqrt{r_0r_1} \sin
\frac{\al^1-\al^0}{2} du\Big)^2\right)^{1/2}
\end{align*}

The distance on $B_{i,\op}$ is the infimum of this expression over
all changes of coordinate for $c^0$. Assuming that $c^0$ and
$c^1$ are originally parametrized with $1/2\pi$ times arc-length so that
$r_0\equiv 1/2\pi, r_1 \equiv 1/2\pi$, this is
\begin{equation}
\label{eq:d.open}
D_{\text{op, diff}}(c^0, c^1) 
= \arccos \sup_\phi \frac{1}{2\pi}\int_0^{2\pi}  \sqrt{\phi_\th} \cos
\frac{\al^1\circ \phi - \al^0}{2} d\th
\end{equation}
and modulo rotations
\begin{multline*}
D_{\text{op, diff, rot}}(c^0, c^1) 
=  \arccos \sup_\phi \left( \Big(\frac{1}{2\pi}\int_0^{2\pi}  
\sqrt{\phi_\th} \cos \frac{\al^1\circ \phi -\al^0}{2} d\th\Big)^2\right.
\\
\left. + \Big(\frac{1}{2\pi}\int_0^{2\pi}  \sqrt{\phi_\th} 
\sin \frac{\al^1\circ \phi -\al^0}{2} d\th\Big)^2\right)^{1/2}.
\end{multline*}
The supremum in both expressions is taken over all increasing bijections 
$\phi\in C^{\infty}([0,2\pi], [0, 2\pi])$.

To shorten these formulae, we will use the following notation. 
Define: 
\begin{align*}
C_-(\phi) &=  
\frac{1}{2\pi}\int_0^{2\pi}  \sqrt{\phi_\th} \cos 
\frac{\al^1\circ \phi -\al^0}{2} d\th
\\
S_-(\phi) &= 
\frac{1}{2\pi}\int_0^{2\pi}  \sqrt{\phi_\th} \sin
\frac{\al^1\circ \phi -\al^0}{2} d\th
\end{align*}
Then we have:  
\begin{align*}
D_{\text{op, diff}} &= \inf_\phi \arccos (C_-(\phi))\quad \text{ and}
\\ 
D_{\text{op, diff, rot}} &=
\inf_{\phi} \arccos\big(\sqrt{(C_-(\phi))^2+(S_-(\phi))^2}\big).
\end{align*}

\subsection{Problems with the existence of geodesics}

These expressions give very explicit descriptions of distance 
and geodesics. We have already noted, however, that, even if 
both $(e^0, f^0)$ and $(e^1, f^1)$ belong to $\onb{S}^{0}$, the 
same property is  not guaranteed at each point of the
geodesic. $e(\al,t) = f(\al, t)=0$ happens for some $t$ whenever
$(e^0(\al), f^0(\al))$ and $(e^1(\al), f^1(\al))$ are collinear 
with opposite orientations. This is not likely to happen for 
geodesics joining nearby points. When it does happen, it is usually
a stable phenomenon: for example, if the geodesic crosses 
$\onb{S}-\onb{S}^{0}$ transversally, as illustrated in figure
\ref{fig:genbifurc}, then this happens for all nearby geodesics too.
Note that this means that the geodesic spray on $\on{Imm}_\op$
is {\it not} surjective. In fact, any geodesic on $\on{Imm}_\op$
comes from a great circle on $\onb{S}^0$ and if it crosses
$\onb{S}-\onb{S}^{0}$, it leaves $\on{Imm}_\op$.

When we pass to the quotient by reparametrizations, another question
arises: does the inf over reparametrizations exist? or equivalently
is there is a {\it horizontal} geodesic joining any two open curves?
In fact, there need not be any such geodesic even if you allow it
to cross $\onb{S}-\onb{S}^{0}$. In general, to obtain a geodesic
minimizing distance between 2 open curves, the curves themselves
must be given parametrizations with zero velocity somewhere, i.e.\
they may need to be lifted to points in $\onb{S}-\onb{S}^{0}$.

\begin{figure}[ht]
\begin{center}
\includegraphics[width=\textwidth]{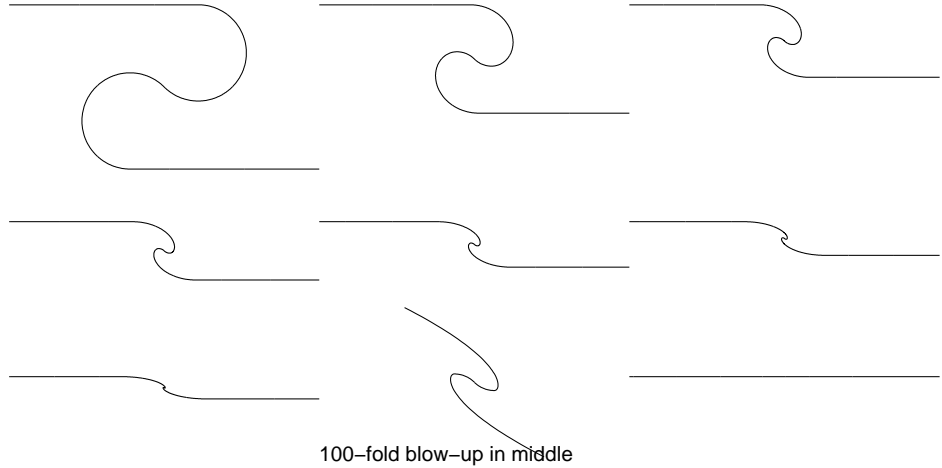}\\ 
\caption{\label{fig:badinf}\footnotesize This is a geodesic 
of open curves running from the curve with the kink at the 
top left to the straight line on the bottom right. A blow up
of the next to last curve is shown to reveal that the kink
never goes away -- it merely shrinks. Thus this geodesic is
not continuous in the $C^1$-topology on $B_\op$. The straight 
line is parametrized so that it stops for a whole interval of 
time when it hits the middle point and thus it is $C^1$-continuous
in $\on{Imm}_\op$.}
\end{center}
\end{figure}

This is best illustrated by the special case in which $c^1$ is the
line segment from $0$ to $1$, namely $e^1+if^1 = 1/\sqrt{\pi}$, 
$\al^1 \equiv 0$. The curve $c^0$ can be arbitrary. Then the 
reparametrization $\phi$ which minimizes distance is the one which
maximizes:
$$ \int_0^{2\pi} \sqrt{\phi_\th} \cos \frac{\al^0}{2} d\th.$$
This variational problem is easy to solve: the optimal $\phi$ is
given by:
$$ \phi(u) = 2\pi \int_0^{u} 
\max\left(\cos \frac{\al^0}{2},0\right)^2 \Bigg / 
\int_0^{2\pi} \max\left(\cos \frac{\al^0}{2},0\right)^2. $$
Note that $\phi$ is {\it not} in general a diffeomorphism: it 
is constant on intervals where $\cos(\al^0/2) \le 0$. Its graph 
is a monotone relation in the sense of section \ref{sec:bigger}.
In fact, it's easy to see that monotone relations enjoy a
certain compactness, so that the inf over reparametrizations
is always achieved by a monotone relation.
Assuming $\al^0$ is represented by a continuous function for which
$-2\pi < \al^0(u) < 2\pi$, the result is that the places on the 
curve $c^0$ where $|\al^0(u)| > \pi$ get squashed to points on the line
segment. The result is that this limit geodesic is not
actually a path in the space $B$ of smooth curves. 
Figure \ref{fig:badinf} illustrates this effect. 

The general problem of maximizing the functional
$$
U(\phi) = \frac{1}{2\pi}\int_0^{2\pi}  \sqrt{\phi_\th} \cos
\frac{\al^1\circ \phi - \al^0}{2} d\th
$$
with respect to increasing functions $\phi$ has been addressed 
in \cite{ty01}. Existence of solutions can be shown in the 
class of monotone relations, or, equivalently functions 
$\phi$ that take the form $\phi(s) = \mu([0,s))$ for some
positive measure $\mu$ on $[0, 2\pi]$ with total mass less 
or equal to $2\pi$ ($\phi_\th$ being replaced by the 
Radon-Nykodim derivative of $\mu$ in the definition of $U$). 
The optimal $\phi$ is a diffeomorphism as soon as 
$\cos ((\al^1(u) - \al^0(v))/{2}) > 0$ whenever $|u-v|$ 
is smaller than a constant (which depends on $\al^0$ and 
$\al^1$).  More details can be found in \cite{ty01}. 

It is easy to show that maximizing $U$ is equivalent to maximizing
$$
U^+(\phi) = \frac{1}{2\pi}\int_0^{2\pi}  \sqrt{\phi_\th} \max\left(\cos
\frac{\al^1\circ \phi - \al^0}{2}, 0\right) d\th
$$
because one can always modify $\phi$ on intervals on which $\cos
((\al^1\circ\phi - \al^0)/{2}) < 0$ to ensure that $\phi_\th d\th =
0$. In \cite{you98}, it is proposed to maximize
$$
\bar U(\phi) = \frac{1}{2\pi}\int_0^{2\pi}  \sqrt{\phi_\th} \left| \cos
\frac{\al^1\circ \phi - \al^0}{2}\right| d\th.
$$
This corresponds to replacing the lift $e(u)+if(u)$ by $\si(u).(e(u)+if(u))$
where $\si(u)=\pm$ for all $u$, but this is beyond the scope of this
paper.

\subsection{Neretin geodesics on $\onb{Gr}(2,V)$}
The integrated path-length distance and explicit geodesics can be
found in any Grassmannian using {\em Jordan Angles} \cite{ner01} as
follows:
If $W_0, W_1 \subset V$ are two 2-dimensional subspaces
the singular value
decomposition of the orthogonal projection $p$ of $W_0$ to $W_1$
gives orthonormal bases $\{e^0, f^0\}$ of $W_0$ and $\{e^1, f^1\}$ of
$W_1$ such that $p(e^0)=\la_e.e^1$, $p(f^0)=\la_f f^1$, $e^0 \perp f^1$
and $f^0 \perp e^1$, where $0 \le \la_f, \la_e \le 1$. Write $\la_e =
\cos(\ps_e), \la_f = \cos(\ps_f)$ then $\ps_e$, $\ps_f$ are the Jordan angles, 
$0\le \ps_e, \ps_f \le \pi/2$. The global metric is given by: 
$$ d(W^0,W^1) = \sqrt{\ps_e^2 + \ps_f^2}$$
and the geodesic by: 
\begin{equation}\label{eq:NeretinGeod}
W(t) = \left\{
{\begin{aligned} e(t)&=\frac{\sin((1-t).\ps_e)e^0 + \sin(t\ps_e)e^1}{\sin\ps_e},
\\
f(t)&=\frac{\sin((1-t).\ps_f)f^0 + \sin(t\ps_f)f^1}{\sin \ps_f}
\end{aligned}}
\right\}.
\end{equation}

We apply this now in order to compute the distance between the 
curves in the two spaces $\on{Imm}_{\ev}/(\text{sim})$ 
and $\on{Imm}_{\od}/(\text{sim})$, as well as in
the unparametrized quotients $B_{i,\ev}/(\text{sim})$ 
and $B_{i,\od}/(\text{sim})$. We write as above  
$\partial_\th c^0 = r_0(\th) e^{i\al^0(\th)}$ and
$\partial_\th c^1 = r_1(\th) e^{i\al^1(\th)}$. 
We put 
\begin{align*}
\bar e^0&=\sqrt{2r_0}\cos\tfrac{\al^0}2 &\quad \bar f^0&=
\sqrt{2r_0}\sin\tfrac{\al^0}2, \\
\bar e^1&=\sqrt{2r_1}\cos\tfrac{\al^1}2 &\quad \bar f^1&=
\sqrt{2r_1}\sin\tfrac{\al^1}2
\end{align*}
thus lifting these curves to 2-planes in the Grassmannian.
The $2\times 2$ matrix of the orthogonal projection 
from the space $\{ \bar e^0, \bar f^0 \}$ to $\{ \bar e^1, \bar f^1 \}$
in these bases is:
$$M(c^0,c^1)=\begin{pmatrix}
\int_{S^1} 2\sqrt{r^0.r^1}. \cos\tfrac{\al^0}2\cos\tfrac{\al^1}2\,d\th & 
  \int_{S^1} 2\sqrt{r^0.r^1}. \cos\tfrac{\al^0}2\sin\tfrac{\al^1}2\,d\th \\ \\
\int_{S^1} 2\sqrt{r^0.r^1} .\sin\tfrac{\al^0}2\cos\tfrac{\al^1}2\,d\th & 
  \int_{S^1} 2\sqrt{r^0.r^1}. \sin\tfrac{\al^0}2\sin\tfrac{\al^1}2\,d\th \end{pmatrix}
$$
It will be convenient to use the notations: 
\begin{align*}
C_\pm &:= \int_{S^1} \sqrt{r^0.r^1}\cos\tfrac{\al^0\pm\al^1}2\,d\th = \tfrac12 \left(M(c^0,c^1)_{11}) \mp M(c^0,c^1)_{22}\right) \\
S_\pm &:= \int_{S^1} \sqrt{r^0.r^1}\sin\tfrac{\al^0\pm\al^1}2\,d\th = \tfrac12 \left(M(c^0,c^1)_{21} \pm M(c^0,c^1)_{12}\right)
\end{align*}

We have to diagonalize this matrix by rotating the curve $c^0$ 
by a constant angle $\be^0$, i.e., the basis $\{\bar e^0,\bar f^0\}$ by the 
angle $\be^0/2$; and similarly the curve $c^1$ by a constant angle 
$\be^1$. So we have to replace $\al^0$ by $\al^0-\be^0$ and 
$\al^1$ by $\al^1-\be^1$ in such a way that 
\begin{align}
0&=\int_{S^1} \sqrt{r^0.r^1} 
  \sin\left(\frac{(\al^0-\be^0) \pm (\al^1-\be^1)}{2}\right) d\th 
  \quad \text{ (for both signs)}
\label{eq:1}
\\&
=S_\pm.\cos\tfrac{\be^0\pm\be^1}2-C_\pm.\sin\tfrac{\be^0\pm\be^1}2 
\notag
\end{align}
Thus
$$ \be_0 \pm \be_1 = 2 \arctan \left( S_\pm / C_\pm \right).$$

In the newly aligned bases, the diagonal elements of $M(c^0,c^1)$ will
be the cosines of the Jordan angles. But even without preliminary 
diagonalization, the following lemma gives you a formula for them:
\begin{lem*} If 
$M=\left( \begin{array}{cc} a & b \\ c & d \end{array}\right)$ 
and $C_\pm = \tfrac12 (a \mp d), S_\pm = \tfrac12(c \pm b)$, 
then the singular values of $M$ are:
$$\sqrt{C_-^2 + S_-^2} \pm \sqrt{C_+^2 + S_+^2}.$$
\end{lem*}

The proof is straightforward. This gives the formula
\begin{align}
D_{\od,\text{rot}}(c^0, c^1)^2 
= &\arccos^2\Big( \sqrt{S_+^2+C_+^2}+ \sqrt{S_-^2+C_-^2}
\Big) \label{eq:d.close,rot} \\ 
+ &\arccos^2\Big( \sqrt{S_-^2+C_-^2}-\sqrt{S_+^2+C_+^2}\Big). \notag
\end{align}
This is the distance in the space 
$\on{Imm}_{\od}(S^1,\mathbb C)$/(\text{transl, rot., scalings}).

\subsection{Horizontal Neretin distances}
If we want the distance in the quotient space 
$B_{i,\od}/(\text{transl, rot., scalings})$ by the group
$\text{Diff}(S^1)$
we have to take the infimum of \thetag{\ref{eq:d.close,rot}} over 
all reparametrizations. To simplify the formulas that follow, we 
can assume that the initial curves $c^0,c^1$ are parametrized by 
arc length so that $r^0 \equiv r^1 \equiv 1/2\pi$. Then consider a
reparametrization $\phi\in\text{Diff}(S^1)$ of one of the two 
curves, say $c^0\o\phi$:
\begin{equation}
\label{eq:d.close,rot,diff}
D_{\text{sim,diff}}(c^0, c^1)^2 
= \inf_{\phi}\left(\arccos^2(\la_e(c^0\o\phi,c^1))
+ \arccos^2(\la_f(c^0\o\phi,c^1))
\right)
\end{equation}
where now  
\begin{align*}
\la_e(c^0\o\phi,c^1) &= \sqrt{S_-^2(\phi)+C_-^2(\phi)}+ \sqrt{S_+^2(\phi)+C_+^2(\phi)}
\\
\la_f(c^0\o\phi,c^1) &= \sqrt{S_-^2(\phi)+C_-^2(\phi)} - \sqrt{S_+^2(\phi)+C_+^2(\phi)}
\\
S_\pm(\phi) &:= \frac{1}{2\pi} \int_{S^1} 
\sqrt{\phi_\th}\sin\tfrac{(\al^0\o\phi)\pm\al^1}2\,d\th,
\\
C_\pm(\phi) &:= \frac{1}{2\pi} \int_{S^1} 
\sqrt{\phi_\th}\cos\tfrac{(\al^0\o\phi)\pm\al^1}2\,d\th.
\end{align*}

To describe the inf, we can use the fact that geodesics on the 
space of curves are the horizontal geodesics in the space of 
immersions. Consider the geodesic $t\mapsto \{e(t),f(t)\}$ in
$\onb{Gr}(2,V)$ described in \thetag{\ref{eq:NeretinGeod}}, for 
\begin{align*}
e^0&=\sqrt{\tfrac{\phi_\th}{\pi}} \cos\tfrac{(\al^0\o\phi)-\be^0}2 &\quad 
e^1&=\tfrac{1}{\sqrt{\pi}} \cos\tfrac{\al^1-\be^1}2, \\
f^0&=\sqrt{\tfrac{\phi_\th}{\pi}} \sin\tfrac{(\al^0\o\phi)-\be^0}2 &\quad 
f^1&=\tfrac{1}{\sqrt{\pi}} \sin\tfrac{\al^1-\be^1}2,
\end{align*}
where the rotations $\be^0$ and $\be^1$ must be computed from 
$c^0\o\phi$ and $c^1$. Note that  
\begin{align*}
e^0_\th &= \frac{\phi_{\th\th}}{2\sqrt{\pi\phi_\th}}
\cos\tfrac{(\al^0\o\phi)-\be^0}2 
- \tfrac1{2\sqrt{\pi}} \phi_\th^{3/2}.(\al^0_\th\o\phi). 
  \sin\tfrac{(\al^0\o\phi)-\be^0}2; 
\\
e^1_\th &= \tfrac{-1}{2\sqrt{\pi}}.\al^1_\th.\sin\tfrac{\al^1-\be^1}2 
\\
f^0_\th &= \frac{\phi_{\th\th}}{2\sqrt{\pi\phi_\th}}
\sin\tfrac{(\al^0\o\phi)-\be^0}2 
+ \tfrac1{2\sqrt{\pi}} \phi_\th^{3/2}.(\al^0_\th\o\phi). 
  \cos\tfrac{(\al^0\o\phi)-\be^0}2; 
\\
f^1_\th &= \tfrac1{2\sqrt{\pi}}.\al^1_\th.\cos\tfrac{\al^1-\be^1}2
\end{align*}
If the Jordan angles are $\ps_e$ and $\ps_f$, then the tangent vector to
the geodesic $t\mapsto W(t)$ at $t=0$ is described by 
$$
e_t(0) = \p_t|_0 e = \frac{\ps_e}{\sin\ps_e}.\left(e^1 - \cos\ps_e.e^0\right),
\;\;
f_t(0) = \p_t|_0 f = \frac{\ps_f}{\sin\ps_f}.\left(f^1 - \cos\ps_f.f^0 \right)
$$
By \ref{horizontalinSt} the geodesic is perpendicular to all
$\on{Diff}(S^1)$-orbits if and only if the sum of Wronskians vanishes:
\begin{align*}
&0=W_\th(e^0,e_t(0))+W_\th(f^0,f_t(0)) =
\\&
=e^0\frac{\ps_e}{\sin\ps_e} \left(e^1_\th - \cos\ps_e e^0_\th \right)
-e^0_\th\frac{\ps_e}{\sin\ps_e} \left(e^1 - \cos\ps_e e^0 \right)
\\&\quad
+f^0\frac{\ps_f}{\sin\ps_f} \left(f^1_\th - \cos\ps_f f^0_\th \right)
-f^0_\th\frac{\ps_f}{\sin\ps_f} \left(f^1 - \cos\ps_f f^0 \right)
\\&
=\frac{\ps_e}{\sin\ps_e} W_\th(e^0,e^1)
+\frac{\ps_f}{\sin\ps_f} W_\th(f^0,f^1)
\\&
=-\frac{1}{\sqrt{\phi_\th}}\Big\{ \phi_{\th\th}
\Big(
 \frac{\ps_e}{\sin\ps_e} \cos\tfrac{(\al^0\o\phi)-\be^0}2 \cos\tfrac{\al^1-\be^1}2
+\frac{\ps_f}{\sin\ps_f} \sin\tfrac{(\al^0\o\phi)-\be^0}2 \sin\tfrac{\al^1-\be^1}2
\Big)
\\
&\quad - \phi_\th \al^1_\th \Big(
\frac{\ps_e}{\sin\ps_e} \cos\tfrac{(\al^0\o\phi)-\be^0}2\sin\tfrac{\al^1-\be^1}2
-\frac{\ps_f}{\sin\ps_f} \sin\tfrac{(\al^0\o\phi)-\be^0}2\cos\tfrac{\al^1-\be^1}2
\Big)
\\
&\quad +\phi_\th^2 (\al^0_\th\o\phi)  \Big(
 \frac{\ps_e}{\sin\ps_e} \sin\tfrac{(\al^0\o\phi)
-\be^0}2\cos\tfrac{\al^1-\be^1}2
 -\frac{\ps_f}{\sin\ps_f} \cos\tfrac{(\al^0\o\phi)-\be^0}2\sin\tfrac{\al^1-\be^1}2
 \Big) \Big\}
\end{align*}
This is an ordinary differential equation for $\phi$ which is coupled to
the (integral) equations for calculating the $\be$'s as functions of $\phi$. 
If it is non-singular (i.e., the coefficient function of 
$\phi_{\th\th}$ does not vanish for any $\th$) then there is a 
solution $\phi$, at least locally. But the non-existence of the
inf described for open curves above will also affect closed curves
and global solutions may actually not exist. However, for closed 
curves that do not double back on themselves too much, as we will see,
geodesics do seem to usually exist.

\subsection{An Example}
\begin{figure}[ht]
\begin{center}
\includegraphics[width=\textwidth]{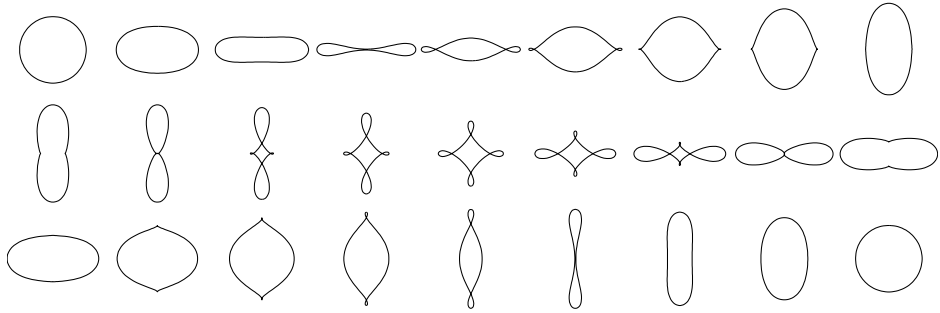}\\ 
\caption{\label{fig:grtcircle}\footnotesize A great circle geodesic on $B_\od$. 
The geodesic begins at the circle at the top left, runs from left to right, 
then to the second row and finally the third. It leaves $B_\od$ twice:
at the top right and bottom left, in both of which the singularity
of figure \ref{fig:genbifurc} occurs in 2 places. The index of the
curve changes from $+1$ to $-3$ in the middle row. See text.}
\end{center}
\end{figure}

Geodesics in the sphere are great circles, which go all the way around the 
sphere and are always closed geodesics. In the case of the Grassmannian, 
using the Jordan angle basis, the geodesic can be continued indefinitely 
using formula \ref{eq:NeretinGeod} above. In fact, it will be a 
{\it closed} geodesic if the Jordan angles $\ps_e, \ps_f$ are 
commensurable. It is interesting to look at an example to see 
what sort of immersed curves arise, for example, at the antipodes to the 
point representing the unit circle. To do this, we take 
$c^0(\th)=e^{i\th}/2\pi$ to be the circle of unit length, giving the 
orthonormal basis $e^0=\cos(\th/2)/\sqrt{\pi}, f^0 = \sin(\th/2)/\sqrt{\pi}$. 
We want $e^1,f^1$ to lie in a direction horizontal with respect to these
and the simplest choice satisfying the Wronskian condition is:
$$ e^1+if^1=\frac{e^{i\th/2}}{\sqrt{\pi}}.\left(\frac{\cos(2\th)}{2} 
- i \sin(2\th)\right).$$
The result is shown in figure \ref{fig:grtcircle}.

\section{Sectional curvature} 

We compute, in this section, the sectional curvature of
$B_{i,\od}/(\text{sim})$ (i.e., translations, rotations,
scaling). We first compute the sectional curvature on the Grassmannian
which is non-negative (but vanishes on many planes) and conclude from 
O'Neill's formula \cite{one66}
that the sectional curvature on $B_i/(\text{sim})$ is non-negative
But since the O'Neill correction term is difficult to compute in this
setting we also do it in a more explicit way, computing first the curvature
on the Stiefel manifold by Gauss' equation, then carrying it over to
$\on{Imm}/(\text{transl})$. Since this is an open subset in a Fr\'echet
space, the O'Neill correction term can be computed more easily on 
$\on{Imm}/(\text{transl})$ and so we finally get a more explicit formula
for the sectional curvature on $B_i/(\text{transl})$. 

\subsection{Sectional curvature on $\onb{Gr}(2,V)$}\label{curvGr}
Let $W\in\onb{Gr}(2,V)$ be a fixed 2-plane which we identify again with 
$\mathbb R^2$. 
Let $\et:V\to V$ be the isomorphism which equals $-1$ on $W$ and $1$ on 
$W^\bot$ satisfying $\et=\et\i$. Then the Grassmanian 
is the symmetric space $O(V)/(O(W)\x O(W^\bot))$ with the involutive
automorphism $\si:O(V)\to O(V)$ given by $\si(U)=\et.U\et$. 
For the Lie algebra in the $V=W\oplus W^\bot$-decomposition we have
$$
\begin{pmatrix} -1 & 0 \\ 0 & 1 \end{pmatrix}
\begin{pmatrix} x & -y^T \\ y & U \end{pmatrix}
\begin{pmatrix} -1 & 0 \\ 0 & 1 \end{pmatrix}
= \begin{pmatrix} x & y^T \\ y & U \end{pmatrix}
$$
Here $x \in L(W,W), y \in L(W,W^\bot)$. The fixed point 
group is $O(V)^\si=O(W)\x O(W^\bot)$. The reductive decomposition 
$\mathfrak g = \mathfrak k + \mathfrak p$ is given by 
$$
\left\{\begin{pmatrix} x & -y^T \\ y & U \end{pmatrix}\right\}
=\left\{\begin{pmatrix} x & 0 \\ 0 & U \end{pmatrix}, 
  x\in \mathfrak s\mathfrak o(2)\right\}
+\left\{\begin{pmatrix} 0 & -y^T \\ y & 0 \end{pmatrix}, 
y\in L(W,W^\bot)\right\}
$$
Let $\pi:O(V)\to O(V)/(O(W)\x O(W^\bot))=\onb{Gr}(2,V)$ be the quotient
projection. Then $T_e\pi:\mathfrak p \to T_o\onb{Gr}$ is an isomorphism, and the
$O(V)$-invariant Riemannian metric on $\onb{Gr}(2,V)$ is given by 
\begin{align*}
G^{\onb{Gr}}_o&(T_e\pi.Y_1,T_e\pi.Y_2)= -\tfrac12\on{tr}(Y_1Y_2)
=-\tfrac12\on{tr}
  \begin{pmatrix} 0 & -y_1^T \\ y_1 & 0 \end{pmatrix}
  \begin{pmatrix} 0 & -y_2^T \\ y_2 & 0 \end{pmatrix}
\\&
=-\tfrac12\on{tr}
  \begin{pmatrix} -y_1^Ty_2 &0  \\ 0 & -y_1y_2^T \end{pmatrix}
= \tfrac12\on{tr}_W(y_1^Ty_2) + \tfrac12\on{tr}_{W^\bot}(y_1y_2^T)
\\&
= \on{tr}_W(y_1^Ty_2)
= \langle y_1(e_1),y_2(e_1) \rangle_{W^\bot}
+ \langle y_1(e_2),y_2(e_2) \rangle_{W^\bot}
\end{align*}
for $Y_1,Y_2\in\mathfrak p$, where $e_1,e_2$ is an orthonormal base of $W$.
By the general theory of symmetric spaces \cite{hel78}, 
the curvature is given by
\begin{align*}
&R^{\onb{Gr}}_o(T_e\pi.Y_1,T_e\pi.Y_2)T_e\pi.Y_1 = T_e\pi.[[Y_1,Y_2],Y_1]
\\&
\left[\begin{pmatrix} 0 & -y_1^T \\ y_1 & 0 \end{pmatrix},
  \begin{pmatrix} 0 & -y_2^T \\ y_2 & 0 \end{pmatrix} \right]
=\begin{pmatrix} -y_1^Ty_2 + y_2^Ty_1 \\ 0 & -y_1y_2^T + y_2y_1^T\end{pmatrix}
\\&
\left[\left[\begin{pmatrix} 0 & -y_1^T \\ y_1 & 0 \end{pmatrix},
  \begin{pmatrix} 0 & -y_2^T \\ y_2 & 0 \end{pmatrix} \right],
  \begin{pmatrix} 0 & -y_1^T \\ y_1 & 0 \end{pmatrix}\right]
=\\&\qquad\qquad
=\left[\begin{pmatrix} -y_1^Ty_2 + y_2^Ty_1 & 0 \\ 
            0 & -y_1y_2^T + y_2y_1^T\end{pmatrix},
  \begin{pmatrix} 0 & -y_1^T \\ y_1 & 0 \end{pmatrix}\right]
\\&\qquad\qquad
= \begin{pmatrix} 0 & 2y_1^Ty_2y_1^T-y_2^Ty_1y_1^T-y_1^Ty_1y_2^T \\ 
             -2y_1y_2^Ty_1+y_2y_1^Ty_1+y_1y_1^Ty_2 & 0 \end{pmatrix}
\end{align*}
For the sectional curvature we have (where we assume that $Y_1,Y_2$ is
orthonormal):
\begin{align*}
k^{\onb{Gr}(2,V)}_{\on{span}(Y_1,Y_2)}&=-B(Y_2,[[Y_1,Y_2],Y_1]) 
= \on{tr}_W(y_2^Ty_2y_1^Ty_1+y_2^Ty_1y_1^Ty_2-2y_2^Ty_1y_2^Ty_1)
\\&
= \tfrac12\on{tr}_W\big((y_2^Ty_1-y_1^Ty_2)^T(y_2^Ty_1-y_1^Ty_2)\big)
\\&\qquad
+ \tfrac12\on{tr}_{W^\bot}\big((y_2y_1^T-y_1y_2^T)^T(y_2y_1^T-y_1y_2^T)\big) 
\\&
= \tfrac12\|y_2^Ty_1-y_1^Ty_2\|^2_{L^2(W,W)}
+ \tfrac12\|y_2y_1^T-y_1y_2^T\|^2_{L^2(W^\bot,W^\bot)}
\ge 0.
\end{align*}
where $L^2$ stands for the space of Hilbert-Schmidt 
operators. Note that there are many orthonormal pairs 
$Y_1,Y_2$ on which sectional curvature vanishes
and that its maximum value 2 is attained when $y_i$ are 
isometries and $y_2=Jy_1$ where $J$ is rotation through 
angle $\pi/2$ in the image plane of $y_1$.

\subsection{Sectional curvature on $\on{Imm/(sim)}$}
The curvature formula can be rewritten by `lowering the indices' which will
make it much easier to express it in terms of the immersion $c$. Fix an
orthonormal basis $e,f$ of $W$ and let $\de e_k = y_k(e), \de f_k = y_k(f)$. 
For $x,y \in W^\bot$, we use the notation $x \wedge y = x\otimes y - y \otimes x 
\in W^\bot \otimes W^\bot.$ Then:
$$ k^{\onb{Gr}(2,V)}_{\on{span}(Y_1,Y_2)}= 
(\langle \de e_1, \de f_2\rangle - \langle \de e_2, \de f_1 \rangle )^2 
+ \tfrac12\|\de e_1 \wedge \de e_2 + \de f_1 \wedge \de f_2\|^2.$$
To check this, note that $y_2^Ty_1-y_1^Ty_2$ is given by a skew-symmetric $2\times 2$
matrix whose off diagonal entry is just 
$\langle \de e_1, \de f_2\rangle - \langle \de e_2, \de f_1 \rangle $ and this identifies 
the first terms in the two formulas for $k$. On the other hand, $y_2y_1^T$
is given by a matrix of rank 2 on the
infinite dimensional space $W^\bot$. 
In view of $W^\bot \otimes W^\bot\subset L(W^\bot,W^\bot)$
it is the 2-tensor 
$\de e_1 \otimes \de e_2 + \de f_1 \otimes \de f_2$. Skew-symmetrizing, we identify
the second terms in the two expressions for $k$.

Going over to the immersion $c$, the tangent vector $\de e_k + i \de f_k$ to $\onb{Gr}$
becomes the tangent vector $h_k = \de c = \int (\de e_k + i\de f_k)(e+if) d\th$
to $\on{Imm/(sim)}$. To express the first term in the curvature, we have:
\begin{prop*}
$\langle \de e_1, \de f_2\rangle - \langle \de e_2, \de f_1 \rangle 
= \int_C \det (D_s h_1, D_s h_2) ds.$
\end{prop*}
\demo{Proof}
$$ D_s(h_k) = \frac{(e+if).(\de e_k + i\de f_k)}{e^2+f^2}$$
hence
$$ \det (D_s h_1, D_s h_2) = \mathrm{Im}( \overline{D_s h_1}, D_s h_2)
= \frac{\mathrm{Im}\left((\de e_1 - i \de f_1)(\de e_2 + i \de f_2)\right)}{e^2+f^2}$$
hence
$$ \int_C \det (D_s h_1, D_s h_2)ds = \int_{S^1} (\de e_1 \de f_2 - \de e_2 \de f_1) d\th.$$
\qed \enddemo
The second term is not quite so compact: because it is a norm on $W^\bot \otimes W^\bot$,
it requires double integrals over $C\times C$, not just a simple integral over $C$.
We use the notation as above $c(\th) = r(\th)e^{i\al(\th)}$. Then we have:
\begin{prop*}
\begin{align*}
&\| \de e_1 \wedge \de e_2 + \de f_1 \wedge \de f_2 \|^2 = 
\mathrm{term 1} + \mathrm{term 2} \\
&\mathrm{term 1} = \iint \limits_{C \times C} \frac{1+\cos(\al(x)-\al(y))}{2} \cdot
\left(\begin{array}{l}\langle D_s h_1(x), D_s h_2(y) \rangle \\
- \langle D_s h_2(x), D_s h_1(y) \rangle \end{array}
\right)^2 ds(x) ds(y) \\
&\mathrm{term 2} = \iint \limits_{C \times C} \frac{1-\cos(\al(x)-\al(y))}{2} \cdot
\left(\begin{array}{l}\det (D_s h_1(x), D_s h_2(y)) \\
- \det (D_s h_2(x), D_s h_1(y)) \end{array}
\right)^2 ds(x) ds(y)
\end{align*}
\end{prop*}
\demo{Proof} Using $r$ and $\al$, we have:
$ \sqrt{r}e^{-i\al/2}D_s h_k = \de e_k + i\de f_k,$
hence:
$$ \sqrt{r(x)r(y)} e^{\frac{i(\al(x)-\al(y)}{2}} \overline{D_s h_1(x)} D_s h_2(y)
= \de e_1(x) \de e_2(y) + \de f_1(x) \de f_2(y) + i(\cdots)
$$
Skew-symmetrizing in the 2 vectors $h_1,h_2$, we get:
\begin{multline*}
\sqrt{r(x)r(y)} \mathrm{Re}\left\{
e^{\frac{i(\al(x)-\al(y))}{2}} \left(\overline{D_s h_1(x)} D_s h_2(y)- 
\overline{D_s h_2(x)} D_s h_1(y) \right)\right\} =
\\ \de e_1(x) \de e_2(y)-\de e_2(x) \de e_1(y) 
+ \de f_1(x) \de f_2(y)-\de f_2(x) \de f_1(y).
\end{multline*}
Squaring and integrating over $S^1\times S^1$, the right hand 
becomes $\|\de e_1 \wedge \de e_2 + \de f_1 \wedge \de f_2\|^2$.
On the left, first write $\mathrm{Re}(
e^{i(\al(x)-\al(y))/2}(\cdots ))$ as the sum of
$\cos((\al(x)-\al(y))/2)\mathrm{Re}(\cdots)$ and
$-\sin((\al(x)-\al(y))/2)\mathrm{Im}(\cdots)$. Then when we
square and integrate, the cross term drops out because it is
odd when $x,y$ are reversed.
\qed \enddemo

We therefore obtain the expression of the curvature in
$\on{Imm/(sim)}$:
\begin{align}
\nonumber
&k^{\on{Imm}/(\text{sim})}_{\mathrm{span}(h_1,h_2)} 
=\Big(\int_C \det (D_s h_1, D_s h_2) ds\Big)^2 +
\\ \label{eq:curv.gr}&\quad
+ \iint \limits_{C \times C} \frac{1+\cos(\al(x)-\al(y))}{2} \cdot
\left(\begin{array}{l}\langle D_s h_1(x), D_s h_2(y) \rangle\\
- \langle D_s h_2(x), D_s h_1(y) \rangle \end{array}
\right)^2 ds(x) ds(y) +
\\\nonumber&\quad
+ \iint \limits_{C \times C} \frac{1-\cos(\al(x)-\al(y))}{2} \cdot
\left(\begin{array}{l}\det (D_s h_1(x), D_s h_2(y)) \\
- \det (D_s h_2(x), D_s h_1(y)) \end{array}
\right)^2 ds(x) ds(y)
\end{align}

A major consequence of the calculation for the curvature on the Grassmannian is:
\begin{thm}
The sectional curvature on $B_i/(\on{sim})$ is $\ge 0$.
\end{thm}

\demo{Proof}
We apply O'Neill's formula \cite{one66} to the Riemannian submersion 
\begin{gather*}
\pi:\onb{Gr}^0\to \onb{Gr}^0/U(V) \cong B_i/\on{Diff}^+(S^1)
\\
k^{\onb{Gr}^0\!\!/U(V)}_{\pi(W)}(X,Y) =
k^{\onb{Gr}^0}_W(X^{\text{hor}},Y^{\text{hor}}) + \tfrac34
\|[X^{\text{hor}},Y^{\text{hor}}]^{\text{ver}}|_{W}\|^2 \ge 0
\end{gather*}
where $X^{\text{hor}}$ is a horizontal vector field projecting to a vector
field $X$ at $\pi(W)$; similarly for $Y$. 
The horizontal and vertical projections exist and are pseudo
differential operators, see \ref{curvB}.
\qed\enddemo

\subsection{Sectional curvature on $\onb{St}(2, V)$}\label{curvStiefel}
The Stiefel manifold is not a symmetric space (as the Grassmannian); it is
a homogeneous Riemannian manifold. This can be used to compute its
sectional curvature. But the following procedure is simpler:

For $(e,f)\in V^2$ we consider the functions 
\begin{equation*}
Q_{1}(e,f)= \frac12 \|e\|^{2}, \quad 
Q_{2}(e,f)= \frac{1}{2} \|f\|^{2}\text{ and }
Q_{3}(e,f)= \frac{1}{\sqrt2}\scp{e}{f}.
\end{equation*}
Then $\onb{St}(2, V)$ is the codimension 3 submanifold 
of $V^2$ defined by the equations $Q_{1}=Q_{2} = 1/2$, 
$Q_{3}=0$.

The metric on $\onb{St}(2, V)$ is induced by the metric 
on $V^{2}$. If ${\xi_1}=(\de e_1, \de f_1)$ and 
${\xi_2}=(\de e_2,\de f_2)$ are
tangent vectors at a point in $V^{2}$, we have
$\scp{{\xi_1}}{{\xi_2}}= 
\scp{\de e_1}{\de e_2} + \scp{\de f_1}{\de f_2}$. 
For a function $\ph$ on $V^{2}$ its gradient 
$\on{grad}\ph$ of $\ph$ (if it exists) is given by 
$\langle \on{grad}\ph(v),\xi \rangle = 
d\ph( v)(\xi) = D_{ v, \xi}\ph$.
The following are the gradients of $Q_i$:
$$
\on{grad} Q_{1} = (e, 0),\quad
\on{grad} Q_{2} = (0, f),\quad
\on{grad} Q_{3} = \frac{1}{\sqrt2} (f, e),
$$
and these form an orthonormal basis of the normal 
bundle $\on{Nor}(\onb{St})$ of $\onb{St}(2, V)$. Let 
${\xi_1},{\xi_2}$ be two normal unit 
vectors tangent to $\onb{St}(2, V)$ at a point
$(e,f)$. Since $V^2$ is flat the sectional 
curvature of $\onb{St}(2, V)$ is given by the Gauss
formula \cite{doc92}: 
\begin{equation*}
k^{\onb{St}(2,V)}_{\mathrm{span}(\xi_1,\xi_2)}
=\langle S({\xi_1},{\xi_1}),
S({\xi_2},{\xi_2})\rangle
-\langle S({\xi_1},{\xi_2}),
S({\xi_1},{\xi_2})\rangle
\end{equation*}
where $S$ denotes the second fundamental form of
$\onb{St}(2, V)$ in $V^{2}$. Moreover, when a
manifold is given as the zeros of functions $F_k$ in
a flat ambient space whose gradients are orthonormal, 
the second fundamental form is given by:
$$ S(X,Y) = \sum_k H_{F_k}(X,Y)\cdot \mathrm{grad} F_k$$
where $H$ is the Hessian of second derivatives. 
Given $\xi_1, \xi_2\in T_{(e,f)}\onb{St}$ with $\xi_i = (\de e_i, \de
f_i)$,
we have:
\begin{gather*}
H_{Q_1}(\xi_1, \xi_2)
  = \langle \de e_1,\de e_2 \rangle,\quad
H_{Q_2}(\xi_1, \xi_2)
  = \langle \de f_1,\de f_2 \rangle, \\
H_{Q_3}(\xi_1, \xi_2) 
  = \tfrac1{\sqrt2}\big(\langle \de e_1,\de f_2 \rangle 
  +\langle \de e_2,\de f_1\rangle\big)
\end{gather*}
so that 
\begin{multline*}
S(\mathbf\xi_1,\mathbf\xi_2) = 
-\langle \de e_1,\de e_2 \rangle \on{grad}Q_1
-\langle \de f_1,\de f_2 \rangle \on{grad}Q_2\\
-\tfrac1{\sqrt2}\big(\langle \de f_1,\de e_2 \rangle 
+\langle \de e_1,\de f_2\rangle\big)
  \on{grad}Q_3.
\end{multline*}
Finally, {\it the sectional curvature 
of $\onb{St}(2,V)$ for  a normal pair of unit
vectors $\mathbf\xi,\mathbf\et$ in $T_{\mathbf f}\onb{St}(2,V)$
is given by:}
\begin{align}
k^{\onb{St}(2,V)}_{\mathrm{span}(\xi_1,\xi_2)}
&= \|\de e_1\|^2 \|\de e_2\|^2 + 
\|\de f_1\|^2 \|\de f_2\|^2 +
2\langle \de e_1,\de f_1 \rangle
\langle\de e_2,\de f_2 \rangle
\notag\\
\label{eq:curvSt} 
&\quad -\langle \de e_1, \de e_2 \rangle^2 - 
\langle \de f_1, \de f_2\rangle^2 
- \tfrac12(\langle \de e_1,\de f_2 \rangle + 
\langle\de f_1,\de e_2 \rangle)^2\notag \\
&= \tfrac12 \| \de e_1 \otimes \de e_2 - \de e_2 \otimes \de e_1 
  + \de f_1 \otimes \de f_2 - \de f_2 \otimes \de f_1 \|^2  \\
&\quad -\tfrac12 (\langle \de e_1, \de f_2 \rangle 
  - \langle \de e_2, \de f_1 \rangle)^2 \notag
\end{align}

Comparing this with the curvature for the Grasmannian,
we see that the O'Neill factor in this case is 
$\tfrac32 (\langle \de e_1,\de f_2 \rangle + 
\langle\de f_1,\de e_2 \rangle)^2$. Moreover, we can write
for the curvature of the isometric $\on{Imm/(transl,scal)}$
\begin{align}
\nonumber
&k^{\on{Imm/(transl,scal)}}_{\mathrm{span}(h_1,h_2)}  
= -\tfrac12 \left(\int_C \det (D_s h_1, D_s h_2) ds \right)^2
+\\
\label{eq:curv.st}
&+\tfrac12
\iint \limits_{C \times C} \frac{1+\cos(\al(x)-\al(y))}{2} \cdot
\left(\begin{array}{l}\langle D_s h_1(x), D_s h_2(y) \rangle - \\
- \langle D_s h_2(x), D_s h_1(y) \rangle \end{array}
\right)^2 ds(x) ds(y)
+\\
\nonumber
&+\tfrac12 
\iint \limits_{C \times C} \frac{1-\cos(\al(x)-\al(y))}{2} \cdot
\left(\begin{array}{l}\det (D_s h_1(x), D_s h_2(y)) - \\
- \det (D_s h_2(x), D_s h_1(y)) \end{array}
\right)^2 ds(x) ds(y) 
\end{align}

\subsection{Sectional curvature on the unscaled Stiefel
manifold}\label{curvunscSt}
Using the basic mapping $\Phi$, the manifold $\on{Imm}/(\text{transl})$ 
can be identified with the unscaled Stiefel manifold which we view as 
the following submanifold of $V^2$. We do not introduce a systematic
notation for it.
\begin{equation}
\label{eq:M}
M = \{(e,f)\in V^2\setminus\{(0,0)\}, 
\|e\|^2 = \|f\|^2 \text{ and }\scp{e}{f} = 0\}
\end{equation}
equipped with the metric
$$
\|(\de e, \de f)\|_{(e, f)}^2 = 2\frac{\|\de e\|^2 +
\|\de f\|^2}{\|e\|^2 + \|f\|^2}.
$$
Consider the diffeomorphism $\Ps: \mathbb R^+ \times \onb{St}(2, V) \to M$ defined by 
$$
\Ps(\ell, (e, f)) =  (\sqrt{\ell}.e,\sqrt{\ell}.f) =: (\bar e,\bar f).
$$
For $\xi (\de e,\de f)\in T_{(e,f)}\onb{St}$ we have
$$
T_{(\ell,e,f)}\Ps.(\la,\xi) =
  \Big(\frac{\la}{2\sqrt{\ell}}e+\sqrt{\ell}.\de e,
  \frac{\la}{2\sqrt{\ell}}f+\sqrt{\ell}.\de f\Big)
  =:(\bar\de e,\bar\de f)
$$
Thus, $\Ps$ is an isometry if $\mathbb R^+ \times \onb{St}(2, V)$ is
equipped with the metric
$$
\|(\la, \xi)\|^2_{\ell, (e, f)} = \frac{\la^2}{2\ell^2}
+ \|\de e\|^2 + \|\de f\|^2
$$
so that $M$ is isometric to the Riemannian product of $\mathbb R^+$ and
$\onb{St}(2, V)$, taking $\|\la\|_\ell = |\la|/(\sqrt2.\ell)$ for the metric on
$\mathbb R^+$. This implies that the curvature tensor on $M$ is the sum of
the tensors on $\mathbb R^+$ (which vanishes) and $\onb{St}(2, V)$. Thus, if 
$\bar\xi_i =  T_{(\ell,f)}\pi.(\la, \xi_i)$,
$i=1,2$ with $(\xi_1, \xi_1)$ orthonormal, 
\begin{align*}
k^M_{\mathrm{span}(\bar \xi_1, \bar \xi_2)} &= 
\frac{-\langle R_M(\bar\xi_1,\bar\xi_2)\bar\xi_1,\bar\xi_2\rangle}
  {\|\bar\xi_1\|^2\|\bar\xi_2\|^2 - \langle\bar\xi_1,\bar\xi_2\rangle^2}
=\frac{-\langle R_{\onb{St}}(\xi_1,\xi_2)\xi_1,\xi_2\rangle}
  {\|\bar\xi_1\|^2\|\bar\xi_2\|^2 - \langle\bar\xi_1,\bar\xi_2\rangle^2}
\\&
=  \frac{k^{\onb{St}}_{\mathrm{span}(\xi_1, \xi_2)}}
{\|\bar\xi_1\|^2\|\bar\xi_2\|^2 - \langle\bar\xi_1,\bar\xi_2\rangle^2}
\end{align*}
Note that we have the relations: 
\begin{eqnarray}\label{eq:xibarxi}
\de \bar e &=& \frac{\la}{2\sqrt{\ell}}e + \sqrt{\ell}.\de e\\
\nonumber
\de \bar f &=& \frac{\la}{2\sqrt{\ell}}f + \sqrt{\ell}.\de f\\
\end{eqnarray}

\subsection{O'Neill's formula} \label{curvB}

For Riemannian submersions, O'Neill formula \cite{one66} states that
the sectional curvature,  in the plane generated by
two horizontal vectors, is given by the curvature computed on the
space ``above'' plus a positive correction term given by $(3/4)$ times
the squared norm of the vertical projection of the Lie bracket of
any horizontal extensions of the
two vectors. We now proceed to the computation of this correction for
the submersion from $\on{Imm}/\mathrm{(sim)}$ to $B_i/\mathrm{(sim)}$.

Because of the simplicity of local charts there, it will be easier to
start from $\on{Imm}/\mathrm{(transl)}$.  Let $c\in \on{Imm}$ with
$\int_{S^1}c\,ds=0$. We first compute the vertical projection of a
vector $h\in T_c\on{Imm}/(\text{transl})$ for the submersion
$\on{Imm}/\mathrm{(transl)} \to B_i/\mathrm{(sim)}$.  
Vectors in the vertical space at $c$ take the form
$$
\tilde h = b v + i \al c + \be c,
$$
each generator corresponding (in this order) to the action of
diffeomorphisms, rotation and scaling ($b$ is a function and $\al,
\be$ are constants).  Denoting $h^\top$ the
vertical projection of $h$, and using the fact that $G_c(h, \tilde h)
= G_c(h^\top, \tilde h)$ for any vertical $\tilde h$, we easily obtain
the  fact that 
$$
h^\top = b v + i \al c + \be c,
$$
with 
\begin{eqnarray*}
L^\top b + \al \ka &=& v\cdot Lh,\\
\avg{b\ka} + \al &=& \avg{D_sh \cdot n},\\ 
\be &=& \avg{D_sh \cdot v},
\end{eqnarray*}
where we have used the following notation: $Lh = -D_s^2 h$, $L^\top b = -D_s^2 b +
\ka^2 b$ and, as before
$$
\avg{F} = \frac{1}{\ell} \int F ds.
$$

 From this, we obtain the fact that $b$ must satisfy
\begin{equation}
\label{eq:b}
L^\top b - \avg{b \ka}\ka = v\cdot Lh - \avg{D_sh\cdot n}\ka.
\end{equation}
The operator $L^\top$ is of order two, unbounded, selfadjoint, and positive
on $\{f\in L^2(S^1,ds): \int f ds=0\}$ thus it is invertible on 
$\{f\in C^\infty(S^1,\mathbb R):\int f\,ds=0\}$ by an index argument as given in 
lemma \cite[4.5]{MM3}. The operator $\tilde L^\top$ in the left-hand
side of \eqref{eq:b} is also invertible under the condition that $c$
is not a circle, with an inverse given by
\begin{equation}
\label{eq:L.inv}
(\tilde L^\top)^{-1}\psi = (L^\top)^{-1}\psi +
\frac{\avg{(L^\top)^{-1}\psi \ka}}{1 - \avg{ \ka (L^\top)^{-1}\ka}}
(L^\top)^{-1} \ka.
\end{equation}
This is well defined unless $\ka \equiv \mathrm{constant}$. Indeed, letting $f =
(L^\top)^{-1}\ka$, we have $-fD_s^2 f + \ka^2 f^2 = \ka f$ which
implies $\avg{\ka f} \geq \avg{\ka^2 f^2}$. By Schwartz inequality
we have $\avg{\ka f} \leq (\avg{\ka^2 f^2})^{1/2}$ which
ensures 
$\avg{\ka f} \leq 1.$ Equality requires $\avg{f D_s^2 f} = 0$ or
$f=\text{constant}$, which in turn implies that $\ka = \text{constant}$ and that
$c$
is a circle. We note for future use that $(\tilde
L^\top)^{-1}\ka = (L^\top)^{-1}\ka$. 

We hereafter assume that $c$ has length 1, is parametrized with
its arc-length divided by $2\pi$, and that it is different from the unit
circle (which is a singular
point in $B_i/(\mathrm{sim})$). We can therefore write
\begin{equation}
\label{eq:h.top}
h^\top = \left((\tilde L^\top)^{-1}\psi(h) v - i\avg{\ka (\tilde
L^\top)^{-1}\psi(h) }c\right)  + i \avg{D_sh\cdot v}c + \avg{D_sh\cdot n} c
\end{equation}
with $\psi(h) = v\cdot Lh - \avg{D_sh\cdot n}\ka$.

The right-hand term in  \eqref{eq:h.top} is the sum of three orthogonal
terms, the last two forming the vertical projection for the submersion 
$\on{Imm}/(\mathrm{transl}) \to \on{Imm}/(\mathrm{sim})$. 
Applying O'Neill's formula two times to this submersion
and to $\on{Imm}/(\mathrm{sim}) \to B_i/(\mathrm{sim})$, 
we see that the correcting term for the sectional curvature on 
$B_i/(\mathrm{sim})$ relative to the curvature
on $\on{Imm}/(\mathrm{sim})$, 
in the direction of the horizontal vectors $h_1$ and $h_2$ is  
$$
\rho(h_1, h_2)_c = \frac{3}{4} 
  \left\| (\tilde L^\top)^{-1}\psi([\bar h_1, \bar h_2]_c) v 
  - i\avg{\ka (\tilde L^\top)^{-1}\psi([\bar h_1, \bar h_2]_c)} c\right\|^2,
$$
$\bar h_1, \bar h_2$ being horizontal extensions of $h_1$ and
$h_2$. From the identity
\begin{eqnarray*}
\|bv - i\avg{\ka b}c\|_c^2 &=& \int |b'v +\ka bn - \avg{\ka b}n|^2 ds \\
&=& \int (b Lb + \ka^2 b^2 n - \avg{\ka b}\ka b) ds \\
&=& \int b(\tilde L^T b) ds
\end{eqnarray*}
we can write
$$
\rho(h_1, h_2)_c = \frac{3}{4} \int \psi([\bar h_1, \bar h_2]_c) 
(\tilde L^\top)^{-1}\psi([\bar h_1, \bar h_2]_c) ds.
$$

We now proceed to the computation of the Lie bracket:
\begin{prop} $\psi([h_1^\bot, h_2^\bot]_c) = W_s(D_sh_1 \cdot
n, D_sh_2\cdot n) -
\avg{\mathrm{det}(D_sh_1, D_sh_2)}\ka$
where $W_s(h,k) = h D_sk - k D_sh$ is the Wronskian with respect to
the arc-length parameter.
\end{prop}
\begin{demo}{Proof}
We take $h_1,h_2\in \{f\in C^\infty(S^1,\mathbb R^2):\int f\,ds=0\}$ which
are horizontal at $c$, consider them as constant vector fields on 
$\on{Imm}/(\text{transl})$ and take, as horizontal extensions,  
their horizontal projections
$\ga \mapsto h_1^\bot(\ga), h_2^{\bot}(\ga)$. 
Then we compute the Lie bracket evaluated
at $\ga$:
\begin{eqnarray*}
[h_1^\bot,h_2^\bot]\Big|_{\ga} &=& D_{c,h_2}h_1^\bot(\ga) -
D_{c,h_1}h_2^\bot(\ga) \\
&=& - D_{c,h_2}h_1^\top(\ga) +
D_{c,h_1}h_2^\top(\ga) 
\end{eqnarray*}
since $h_i^\top + h_i^\bot = h_i$ is constant, for $i=1,2$. We have
\begin{eqnarray*}
h_1^\top(\ga) &=& \left((\tilde L^\top_\ga)^{-1}\psi_\ga(h_1) v_\ga - i\avg{\ka_\ga (\tilde
L^\top_\ga)^{-1}\psi_\ga }_\ga \ga \right) \\
&+& i \avg{D_{s_\ga} h_1 \cdot n_\ga}_\ga\ga
+ \avg{D_{s_\ga} h_1 \cdot v_\ga}_\ga \ga
\end{eqnarray*}
with $\psi_\ga(h_1) = v_\ga \cdot L_\ga h_1 - \avg{D_{s_\ga} h_1 \cdot
v_\ga }_\ga \ka_\ga$. We have added subscripts $\ga$ to quantities that
depend on the curve, with $D_{s_\ga}$ holding for the derivative with
respect to the $\ga$ arc-length (we still use no subscript for $\ga =
c$). Note that $\avg{D_{s_\ga} h_1 \cdot n_\ga}_\ga = \ell_\ga \avg{D_s h_1 \cdot
n_\ga}$ and $\avg{D_{s_\ga} h_1 \cdot v_\ga}_\ga = \ell_\ga \avg{D_s h_1 \cdot
v_\ga}$ which is a first simplification. Also, since we assume that
$h_1$ is horizontal at $c$, we have $\avg{D_s h_1 \cdot n} = \avg{D_s
h_1 \cdot v} = 0$ and $v\cdot Lh_1 = 0$, which imply 
$\psi(h_1)=0$. 

We therefore have (to simplify, we temporarily
use the notation $f' = D_sf$)
\begin{eqnarray}
\label{eq:dh.1}
D_{c, h_2} h_1^\top &=& \left((\tilde L^\top)^{-1} D_{c,h_2} \psi_\ga(h_1) v - i\avg{\ka (\tilde
L^\top)^{-1} D_{c,h_2} \psi_\ga(h_1) } c\right) \\
\nonumber 
& +& i \avg{h'_1
\cdot D_{c,h_2}n_\ga} c
+ \avg{h'_1 \cdot D_{c, h_2} v_\ga} c.
\end{eqnarray}

Since we have $D_{c,h_2} v_\ga = (h'_2\cdot n) n $ and $D_{c,h_2} n_\ga =
- (h'_2\cdot n) v$ we immediately obtain the expression of the last two
terms in \eqref{eq:dh.1}, which are
\begin{equation}
\label{eq:dh.2}
- i \avg{(h'_1\cdot v) (h'_2\cdot n)} c
+ \avg{(h'_1\cdot n) (h'_2\cdot n)} c.
\end{equation}

We now focus on the variation of $\phi_\ga$. We need to compute
\begin{eqnarray*}
D_{c,h_2} \psi_\ga(h_1) &=& D_{c,h_2} (v_\ga\cdot L_\ga h_1) -
\avg{(h'_1\cdot v) (h'_2\cdot n)}\ka.
\end{eqnarray*}
If $h$ is a constant vector field, we have 
$$
D_{s_\ga} h = h' \|\ga'\|^{-1}
$$
and 
$$
L_\ga h = - (h' \|\ga'\|^{-1})' \|\ga'\|^{-1}. 
$$
This implies
\begin{eqnarray*}
D_{c, h_2} L_\ga h_1 
&= - h''_1 D_{c,h_2}\|\ga'\|^{-1} - (h'_1 D_{c,h_2}\|D_s\ga\|^{-1})' \\
&=  2h''_1 (h'_2\cdot v) +  h'_1 (h'_2 \cdot v)' .
\end{eqnarray*}
Therefore
\begin{eqnarray*}
D_{c, h_2} (L_\ga h_1 \cdot v_\ga) 
&=   (h'_1\cdot v) (h'_2 \cdot v)' -  (h''_1 \cdot n) (h'_2 \cdot
n).
\end{eqnarray*}
Using
\begin{eqnarray*}
h''_1 &= ((h'_1\cdot v) v + (h'_1\cdot n)n)'\\
&= ((h'_1\cdot v)'- \ka (h'_1\cdot n)) v + ((h'_1\cdot n)'
+ \ka (h'_1\cdot v))n
\end{eqnarray*}
and the fact that $h''_1 \cdot v = h''_2 \cdot v = 0$, we can write
\begin{eqnarray*}
D_{c, h_2} (L_\ga h_1 \cdot v_\ga) 
&=   -  (h'_2 \cdot
n) (h'_1\cdot n)' 
\end{eqnarray*}
which yields
$$
D_{c,h_2} \psi_\ga(h_1) = -  (h'_2 \cdot
n) (h'_1\cdot n)' - \avg{(h'_1\cdot v) (h'_2\cdot n)}\ka.
$$
By symmetry
$$
D_{c,h_2} \psi_\ga(h_1) - D_{c,h_1} \psi_\ga(h_2) = W_s(h'_1 \cdot
n, h'_2\cdot n) - \avg{\text{det}(h'_1, h'_2)}\ka,
$$
where 
$W_s(\varphi_1, \varphi_2) =\varphi _{1}\varphi _{2}^{\prime }
-\varphi_{1}^{\prime }\varphi_{2}.$

Combining this with \eqref{eq:dh.2}, we get  
\begin{eqnarray*}
[h^\bot_1, h_2^\bot]_c &=& (\tilde L^\top)^{-1} (W_s(h'_1 \cdot
n, h'_2\cdot n)  -
\avg{\mathrm{det}(h'_1, h'_2)}\ka) v\\
&-& i\avg{\ka (\tilde
L^\top)^{-1}  (W_s(h'_1 \cdot
n, h'_2\cdot n) -
\avg{\mathrm{det}(h'_1, h'_2)}\ka)} c - i\avg{\mathrm{det}(h'_1, h'_2)} c
\end{eqnarray*}
so that
\begin{eqnarray*}
\psi([h_1^\bot, h_2^\bot]_c) &=& W_s(h'_1 \cdot
n, k'_2\cdot n) -
\avg{\mathrm{det}(h'_1, h'_2)}\ka. 
\end{eqnarray*}
(We have used the fact that $\psi(b v + i\al c) = \tilde L^\top
b$.) 
\qed \end{demo}
We therefore obtain the formula 
\begin{multline}
\label{eq:one.sim}
\rho(h_1, h_2)_c = \frac{3}{4} \int \Big(W_s(D_sh_1 \cdot
n, D_sh_2\cdot n) -
\avg{\mathrm{det}(D_sh_1, D_s h_2)}\ka\Big)\cdot \\ 
\cdot (\tilde L^\top)^{-1} \Big(W_s(D_sh_1 \cdot n, D_sh_2\cdot n) -
\avg{\mathrm{det}(D_sh_1, D_sh_2)}\ka\Big) ds
\end{multline}
with $(\tilde L^\top)^{-1}$ given by \eqref{eq:L.inv}. Finally,
assuming that $h_1$ and $h_2$ are orthogonal,
\begin{align*}
k_{\mathrm{span}(h_1,h_2)}^{\on{B_i/(\text{sim})}} &= 
k_{\mathrm{span}(h_1,h_2)}^{\on{Imm/(\text{sim})}} + \rho(h_1,h_2)_c
\end{align*}
where $k_{\mathrm{span}(h_1,h_2)}^{\on{Imm/(\text{sim})}}$ is given in 
\eqref{eq:curv.gr}.

\bigskip

A similar (and simpler) computation provides the correcting term for the
space $B_i/(\mathrm{transl, scale})$. In this case, the rotation part of the
vertical space disappears, and the two remaining components
(parametrization and scale) are orthogonal. The result is
\begin{equation}
\label{eq:one.scale}
\rho(h_1, h_2)_c = \frac{3}{4} \int W_s(D_sh_1 \cdot
n, D_sh_2\cdot n) (L^\top)^{-1} W_s(D_sh_1 \cdot
n, D_sh_2\cdot n)  ds.
\end{equation}

\subsection{An upper bound for $k^{\on{B_i/(\text{sim})}}_{\mathrm{span}(h_1,h_2)}$}

Here we derive an explicit upper bound for 
$k^{B_i/(\text{sim})}_{\mathrm{span}(h_1,h_2)}$ 
at a fixed curve $c\in \on{B_i/(sim)}$
and a fixed tangent vector $h_2$. This will show that geodesics (such
as the one in the $h_1$ direction)  have
at least a small interval before they meet another geodesic. The size
of this interval can be controlled, as we will see, by an upper bound
that involves the supremum norm of the first two derivatives of $h_1$.

 We assume that
$c$ has length $2\pi $. Since $\on{Imm}/(\text{sim})$ is isometric to
$\onb{Gr}(2,V)$, its sectional curvature is not larger than 2 as already
remarked.
We estimate the terms in $\rho
(h_1,h_2)_{c}=\frac{3}{4}\overline{\left\langle \psi
(h_1,h_2)(\widetilde{L}^{\top })^{-1}\psi (h_1,h_2)\right\rangle }$
where
$\psi (h_1,h_2)=W_s(D_sh_1\cdot n,D_sh_2\cdot
n)-\overline{\left\langle \det (D_s h_1, D_sh_2)\right\rangle }\kappa $. 
For a fixed $h_2$, $\psi (h_1,h_2)$ is function
of $h_1$ belonging to $H^{-1}(c)$. We estimate $\|\psi (h_1,h_2)\|_{c,-1}$
and then $\rho (h_1,h_2)_{c}$ by estimating the norm of the
operator $(L^{\top })^{-1}$ which maps $H^{-1}(c)$to $H^{1}(c)$. 

If $f\in H^{0}(c)$, $\|f\|_{c,-1}\leq \|f\|_{c,0}$ and $\|f^{\prime
}\|_{c,-1}\leq \|f\|_{c,0}$.  Therefore,
\begin{multline*}
\|W_s(D_sh_1\cdot n,D_s h_2\cdot n)\|_{-1}
=\|(D_sh_1\cdot n)D_s(D_sh_2\cdot n)-D_s(D_sh_1\cdot n)(D_sh_2 \cdot n)\|_{-1} \\
\leq \left(
\|D_sh_2\cdot n\|_{c,\infty }+\|D_s(D_sh_2\cdot n)\|_{c,\infty }\right) \|D_sh_1\cdot
n\|_{c,0}.
\end{multline*}
Since $h_1$ has norm 1, $\|D_sh_1\cdot n\|_{c,0}$ and $\|D_sh_1 \cdot v\|_{c,0}$ are $\leq
\sqrt{2\pi }$.
\begin{multline*}
\overline{\left\langle \det
(D_s h_1, D_sh_2)\right\rangle } \leq
\overline{\left\langle |D_s h_1\cdot n|\cdot |D_s h_2\cdot
v|\right\rangle }+\overline{\left\langle |D_sh_1\cdot v|\cdot
|D_sh_2\cdot n|\right\rangle } \\
\leq \frac{1}{2\pi }\left(
\|D_s h_1\cdot n\|_{c,0}\cdot \|D_sh_2\cdot
v\|_{c,0}+\|D_sh_1\cdot v\|_{c,0}\cdot \|D_s h_2\cdot
n\|_{c,0}\right) \leq 2.
\end{multline*}
This results in 
\[\|\psi
(h_1,h_2)\|_{c,-1}\leq \sqrt{2\pi }\left( \|D_s h_2\cdot
n\|_{c,\infty }+\|D_s(D_s h_2\cdot n)\|_{c,\infty
}+2\sqrt{\overline{\left\langle \kappa ^{2}\right\rangle }}\right).
\]
Now 
$$
\overline{\left\langle \psi(\widetilde{L}^{\top })^{-1}\psi \right\rangle}
=\overline{\left\langle \psi (L^{\top })^{-1}\psi \right\rangle}
+\frac{\overline{\left\langle \psi (L^{\top })^{-1}\kappa\right\rangle }^{2}}
{1-\overline{\left\langle \kappa (L^{\top})^{-1}\kappa \right\rangle }}
\leq \frac{\overline{\left\langle
\psi (L^{\top })^{-1}\psi \right\rangle }}{1-\overline{\left\langle
\kappa (L^{\top })^{-1}\kappa \right\rangle }}
$$
since
$\overline{\left\langle \psi (L^{\top })^{-1}\kappa \right\rangle
}^{2} \leq \overline{\left\langle \psi (L^{\top })^{-1}\psi
\right\rangle }\cdot \overline{ \left\langle \kappa (L^{\top
})^{-1}\kappa \right\rangle }$.

\begin{prop}
If $\psi \epsilon H^{-1}(c)$ then
$$ 
\overline{\left\langle \psi(L^{\top })^{-1}\psi \right\rangle }\leq
\frac{1}{2\pi }\left( 1+3\|1-\kappa^{2}\|_{c,\infty }\right) \|\psi
\|_{c,-1}^{2}. 
$$
\end{prop}

\begin{demo}{Proof}
Let $L_{o}=-D_{s}^{2}+1$. Let $L^{\top }f=L_{o}f_{o}=\psi $. Then,
$f,f_{o}\epsilon H^{1}(c)$ and $\|f_{o}\|_{c,1}=\|\psi
\|_{c,-1}$. Let $g=f-f_{o}$ so that $L^{\top }g=(1-\kappa
^{2})f_{o}$. The eigenvalues of $L^{\top }$ are positive and bounded
from below by $1/2$, see \cite{bl04}.  Therefore, $\|g\|_{c,0}^{2}\leq
2(g,L^{\top }g)$ where $(g,L^{\top }g)=\int gL^{\top}gds $. We also
have $\|g^{\prime }\|_{c,0}^{2}\leq (g,L^{\top }g)$ Hence,
$$
\|g\|_{c,1}^{2}\leq 3(g,L^{\top }g)\leq 3\int (1-\kappa
^{2})gf_{o}ds\leq3\|1-\kappa ^{2}\|_{c,\infty
}\|g\|_{c,1}\|f_{o}\|_{c,1}. 
$$
Therefore, 
$\|g\|_{c,1}\leq 3\|1-\kappa^{2}\|_{c,\infty }\|\psi \|_{c,1}$ 
and 
$\|f\|_{c,1}\leq \left(
1+3\|1-\kappa^{2}\|_{c,\infty }\right) \|\psi\|_{c,-1}$.
Finally,
$$
\overline{\left\langle \psi (L^{\top })^{-1}\psi\right\rangle } 
\leq \frac{1}{2\pi }\|\psi \|_{c,-1}\cdot\|(L^{\top})^{-1}\psi \|_{c,1} 
\leq \frac{1}{2\pi }\left( 1+3\|1-\kappa^{2}\|_{c,\infty }\right) \|\psi\|_{c,-1}^{2} 
\qed$$
\end{demo}

Putting all the estimates together we get, for orthonormal $h_1,h_2$ as
always,
\begin{align}
0&\leq k^{\on{B_{i}/(\mathrm{sim})}}_{\mathrm{span}(h_1,h_2)} \leq
\notag\\&
\leq 2+\frac{3\left( 1+3\|1-\kappa
^{2}\|_{c,\infty }\right) \left( \| D_sh_2\cdot n\|_{c,\infty
}+\|(D_sh_2\cdot n)^{\prime }\|_{c,\infty
}+2\sqrt{\overline{\left\langle \kappa^{2}\right\rangle }}\right)
^{2}}{4\left( 1-\overline{\left\langle \kappa(L^{\top })^{-1}\kappa
\right\rangle }\right) }
\end{align}

\section{Numerical procedure and experiments}
The distance $D_{\mathrm{op,dif}}$ given in \ref{greatcircles}  can be
computed in a very short time by
dynamic programming, using a slightly modified procedure from the one
described in \cite{ty00}. Here is a sketch of how it works.

Let $F(\al^0, \al^1) = \max(0, \cos((\al^0-\al^1)/2))$, 
and assume that the curves are discretized
over intervals $[\th^i(k), \th^i(k+1))$, $k=0, n^i-1$, $i=0,1$, so
that the angles have constant values, $\al^i(k)$ on these intervals.
The problem is then equivalent to maximizing
$$
\sum_{k,l} F_{kl} \int_{\max(\th^0(k),
\phi^{-1}(\th^1(l)))}^{\min(\th^0(k+1), \phi^{-1}(\th^1(l+1)))}
\sqrt{\phi_\th}\,d\th
$$
with $F_{kl} = F(\al^0(k), \al^1(l))$. Because the integral of $\sqrt{\phi_\th}$
is maximal for linear $\phi$, we must in fact maximize
\begin{multline*}
\sum_{k,l} F_{kl} \sqrt{\left(\max(\th^0(k), \tilde \th^1(l))) -
\min(\th^0(k+1), \tilde \th^1(l+1)))\right)^+} \\
\sqrt{(\max\left(\tilde \th^0(k), \th^1(l))) -
\min(\tilde \th^0(k+1), \th^1(l+1)))\right)^+}
\end{multline*}
with the notation $\tilde \th^0(k) = \phi(\th^0(k))$ and $\tilde
\th^1(l) = \phi^{-1}(\th^1(l))$. 
The method now essentially implements a coupled linear programming
procedure over the values of $\tilde \th_0$ and $\tilde \th_1$.  See
\cite{you98,ty00} for more details. 
This procedure is very fast, and one still obtains an efficient
procedure by combining it with an exhaustive search for
an optimal rotation.

For closed curves, we can furthermore optimize the result with
respect to  the offset $\phi(0) \in S^1$, for the
diffeomorphism. Doing so provides the value of
$$
D'_{\mathrm{op,diff, rot}}(c^0, c^1) = \inf_\phi \arccos\sqrt{(C_-(\phi))^2 +
(S_-(\phi))^2}.
$$
where the notation $D'$ is here to remember that the minimization is
over $\phi\in C^{\infty,+}(S^1)$ and not $C^{\infty, +}([0, 2\pi])$.

This combination of the almost instantaneous dynamic programming
method and of an exhaustive search over two parameters provides a
feasible elastic matching method for closed curves. But this does not provide
the geodesic distance over ${B}_i/\on{(sim)}$, since we worked with
great circles instead of the Neretin geodesics. There are two
consequences for this: first, the obtained distance is only a lower
bound of the distance on ${B}_i/\on{(sim)}$, and second, since the
closedness constraint is not included, the curves generally become
open during the evolution (as shown in the experiments).

\begin{figure}[ht]
 \begin{center}
\includegraphics[width=0.8\textwidth]{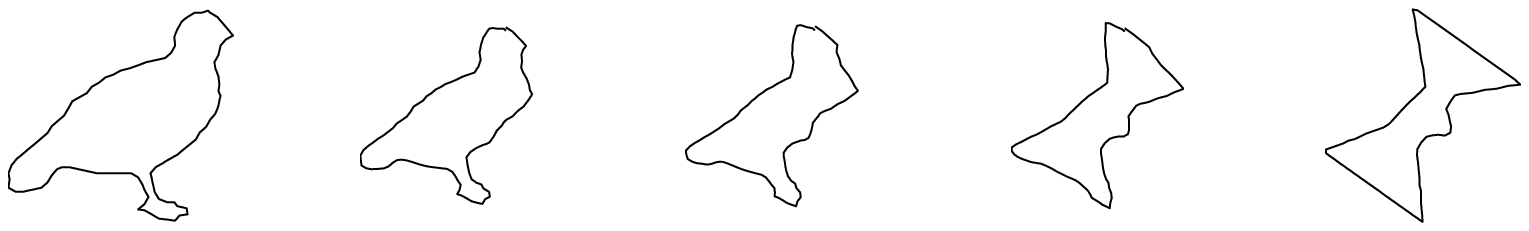}\\ 
\includegraphics[width=0.8\textwidth]{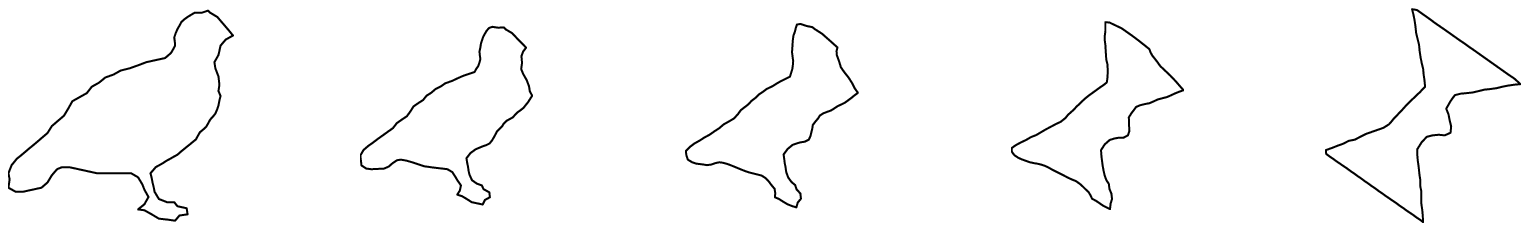} 
\caption{\label{fig:bird02.arb2} Curve evolution with and without the
closedness constraint. Lower and upper bounds for the geodesic distance:
0.443 and 0.444}
\end{center}
\end{figure}

\begin{figure}[ht]
 \begin{center}
\includegraphics[width=0.8\textwidth]{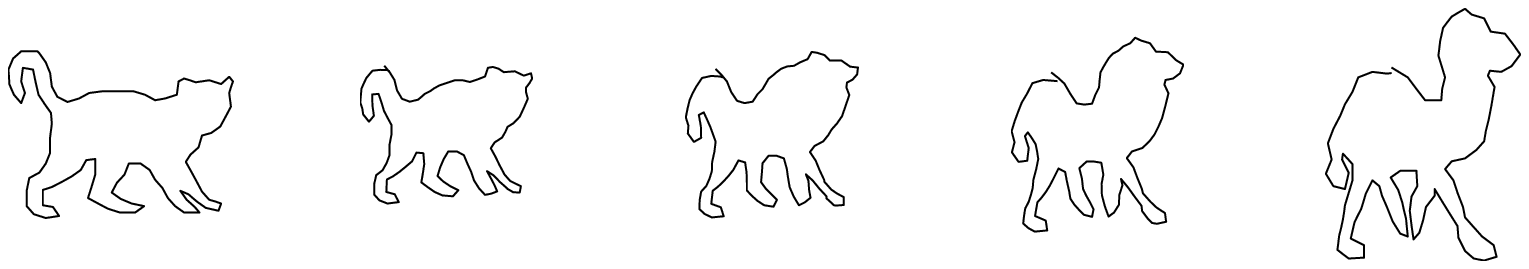}\\ 
\includegraphics[width=0.8\textwidth]{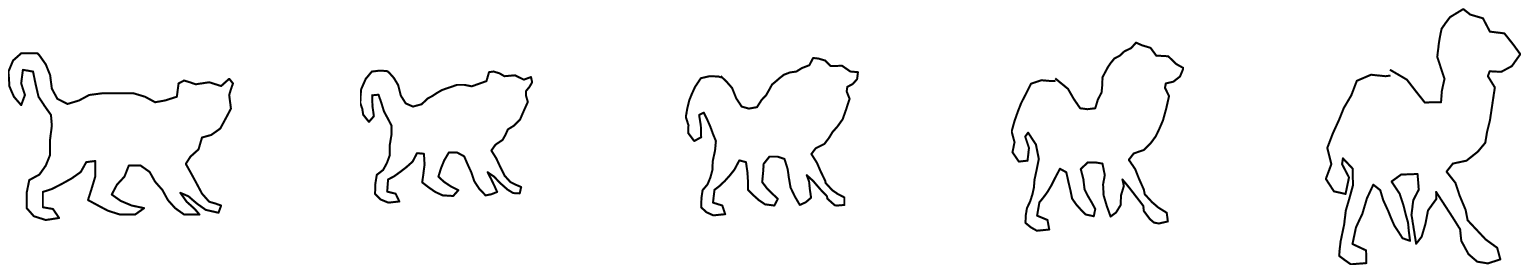} 
\caption{\label{fig:cat02.camel17} Curve evolution with and without the
closedness constraint. Lower and upper bounds for the geodesic distance:
0.462 and 0.464}
\end{center}
\end{figure}

\begin{figure}[ht]
 \begin{center}
\includegraphics[width=0.8\textwidth]{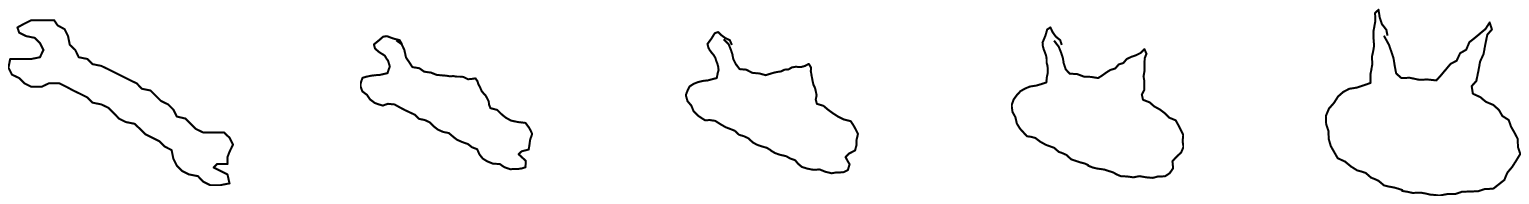}\\ 
\includegraphics[width=0.8\textwidth]{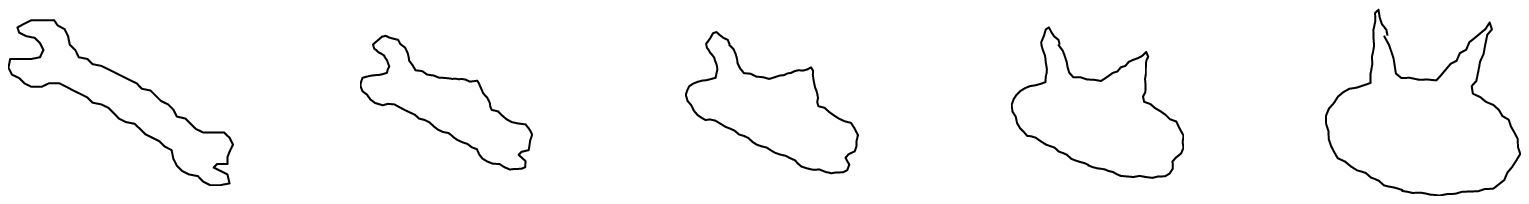} 
\caption{\label{fig:tool07.textbox34} Curve evolution with and without the
closedness constraint. Lower and upper bounds for the geodesic distance:
0.433 and 0.439}
\end{center}
\end{figure}

\begin{figure}[ht]
 \begin{center}
\includegraphics[width=0.8\textwidth]{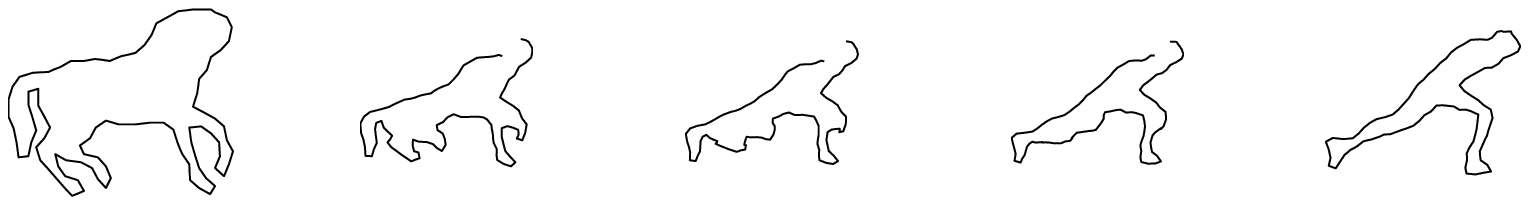}\\ 
\includegraphics[width=0.8\textwidth]{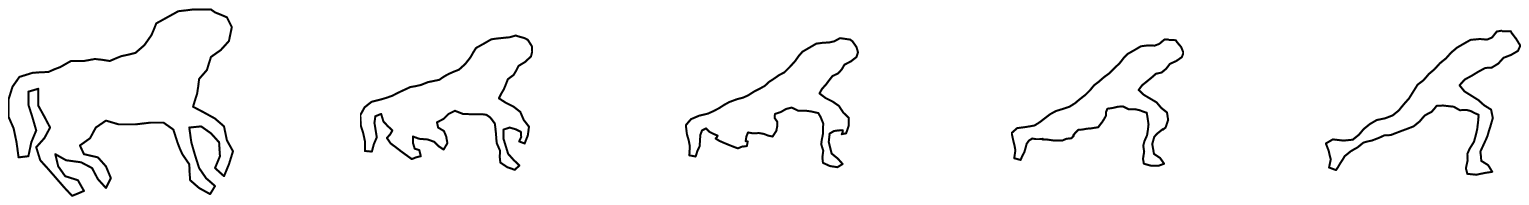} 
\caption{\label{fig:horse05.stef07} Curve evolution with and without the
closedness constraint. Lower and upper bounds for the geodesic distance:
0.498 and 0.532}
\end{center}
\end{figure}

\begin{figure}[ht]
 \begin{center}
\includegraphics[width=0.8\textwidth]{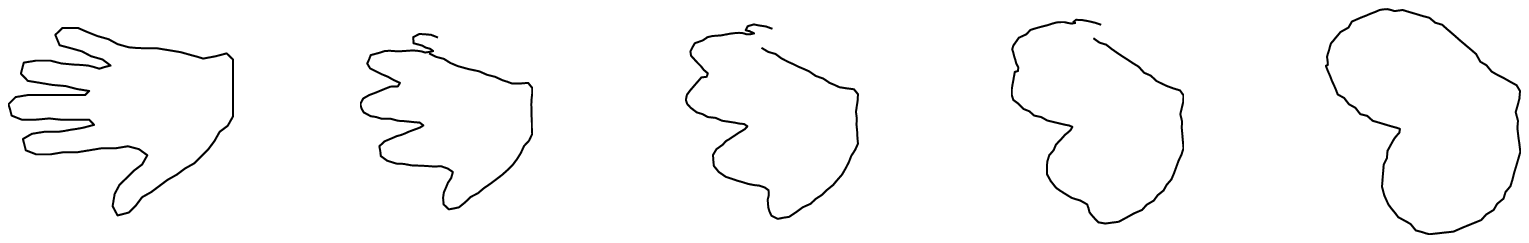}\\ 
\includegraphics[width=0.8\textwidth]{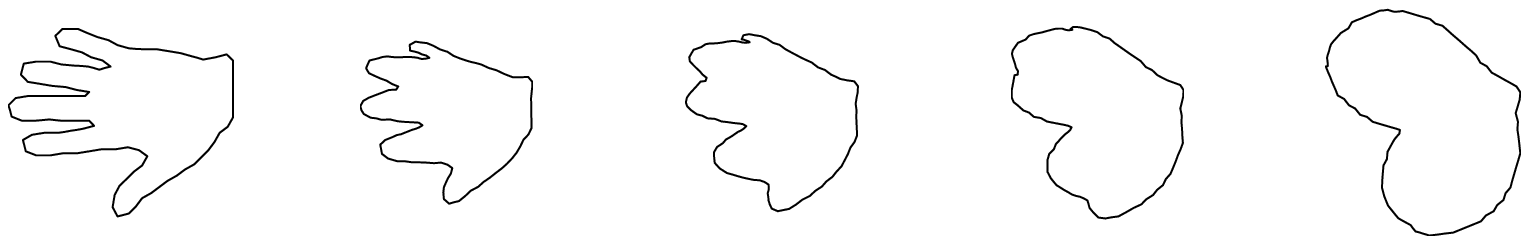} 
\caption{\label{fig:hand90.Heart02} Curve evolution with and without the
closedness constraint. Lower and upper bounds for the geodesic distance:
0.513 and 0.528}
\end{center}
\end{figure}

However, the optimal diffeomorphism which has been obtained by this
approach can be used to reparametrize the curve $c^0$, and we can
compute the geodesic between $c^0\circ \phi^*$ and $c^1$ in
$\on{Imm/(sim)}$ using Neretin geodesics, which forms, this time, an
evolution of closed curves. Its geodesic length now obviously provides
an upper-bound for the geodesic distance on $B_i/\on{(sim)}$. The
numerical results that are presented in figures \ref{fig:bird02.arb2}
to \ref{fig:hand90.Heart02} compare the great circles and Neretin
geodesics obtained using this method. Quite surprisingly,  the
differences between the lower  and upper bounds in these examples are quite
small. 

\newpage

\section{Appendix: The geodesic equation on $G^{\mathrm{imm,scal,1,\infty}}$}

\subsection{The geodesic equation}\label{geodesicImm}
We use the method of \cite{MM3} for the space $\{c\in \on{Imm}_{c}:
c(1)=0\}$ which is an open subset in a Fr\'echet space, with tangent space
$\{h\in C^\infty(S^1,\mathbb C): f(1)=0\}$. 
We shall use the
following conventions and results from \cite{MM3}:
\begin{align*}
&D_s=D_{s,c}=\frac{\p_\th}{|c_\th|},\quad ds = |c_\th|\,d\th,\quad
D_{c,m}\ell(c) = \int \langle D_s m,v_c \rangle ds 
  = -\int \ka_c\langle m,n_c \rangle\,ds,
\\
&D_{c,m}(D_s) = -\langle D_s m, v_c \rangle D_s,\qquad
D_{c,m}(ds) = \langle D_s m,v_c \rangle ds.
\end{align*}
Then the derivative of the metric at $c$ in direction $m$ is:
\begin{align*}
D_{c,m} G_c(h,k) &= 
\frac1{\ell_c^2}\int\ka_c\langle m,n_c \rangle\,ds\cdot 
  \int\langle -D_s^2h,k \rangle\,ds
\\&\quad
+\frac1{\ell_c}\int\langle D_s m,v_c \rangle\langle D_s^2h,k \rangle\,ds
+\frac1{\ell_c}\int\Big\langle D_s\Big(\langle D_s m,v_c\rangle D_s
h\Big),k\Big\rangle\,ds
\\&\quad
-\frac1{\ell_c}\int\langle D_s^2h,k \rangle\langle D_s m,v_c \rangle\,ds
\\&
= \frac1{\ell_c^2}\int\ka_c\langle m,n_c \rangle\,ds\cdot 
  \int\langle -D_s^2h,k \rangle\,ds
\\&\quad
-\frac1{\ell_c}\int\Big\langle -D_s^2m,
  D_s\i\Big(\langle D_s h,D_s k\rangle v_c\Big)\Big\rangle\,ds
\end{align*}
According to \cite[2.1]{MM3} we should rewrite this as
\begin{align*}
D_{(c,m)}G_c(h,k) &=
G_c(K^n_c(m,h), k)
=G_c\big(m, H^n_c(h,k)\big),
\end{align*}
and thus we find the two versions $K$ and $H$ of the $G$-gradient of
$c\mapsto G_c(h,k)$: 
\begin{align*}
K_c(m,h) &= \frac1{\ell_c}\int\ka_c\langle m,n_c \rangle\,ds\cdot h
- D_s\i\big(\langle D_s m,v_c\rangle D_s h\big)
\\
H_c(h,k) &=-\frac1{\ell_c} D_s^{-2}\big(\ka_c n_c\big)\cdot 
  \int\langle -D_s^2h,k \rangle\,ds
-D_s\i\Big(\langle D_s h,D_s k\rangle v_c\Big)
\end{align*}
which gives us the geodesic equation by \cite[2.4]{MM3}:
\begin{align}
c_{tt}&= \tfrac12 H^n_c(c_t,c_t) - K^n_c(c_t,c_t).
\notag\\& 
= -\tfrac12 D_s^{-2}\big(\ka_c n_c\big)\|c_t\|_{G_c}^2 
-\tfrac12 D_s\i\big(|D_s c_t|^2 v_c\big)
\notag\\&\quad
- \frac1{\ell_c}\int\ka_c\langle c_t,n_c \rangle\,ds\cdot c_t
- D_s\i\big(\langle D_s c_t,v_c\rangle D_s c_t\big)
\label{eq:geod}\end{align}

\begin{thm}\label{existenceG}
For each $k\ge 3/2$ the geodesic equation derived in \ref{geodesicImm} 
has unique local solutions in the Sobolev space of
$H^{k}$-immersions. The solutions depend $C^\infty$ on $t$ and on the initial
conditions $c(0,\;.\;)$ and $c_t(0,\;.\;)$.  
The domain of existence (in $t$) is uniform in $k$ and thus this
also holds in $\on{Imm}_*:=\{c\in\on{Imm}(S^1,\mathbb R^2):c(1)=0\}$.
\end{thm}

\begin{demo}{Proof} The proof is very similar to the one of
\cite[4.3]{MM3}. 
We denote by $*$ any space of based loops $(c(1)=0)$.
We consider the geodesic equation as the flow equation of a smooth
($C^\infty$) vector field on the $H^2$-open set $U^k\x H^k_*(S^1,\mathbb R^2)$ 
in the Sobolev space $H^k_*(S^1,\mathbb R^2)\x H^k_*(S^1,\mathbb R^2)$ where
$U^k=\{c\in H^k_*:|c_\th|>0\}\subset H^k$ is $H^2$-open. 
To see that this works we will use the following facts: 
By the Sobolev inequality we have a
bounded linear embedding $H^{k}_*(S^1,\mathbb R^2)\subset C^m_*(S^1,\mathbb R^2)$ if
$k>m+\frac12$. The Sobolev space $H^k_*(S^1,\mathbb R)$ is a Banach algebra
under pointwise multiplication if $k>\frac12$. For any fixed smooth mapping
$f$ the mapping $u\mapsto f\o u$ is smooth $H^k_*\to H^k$ if $k>0$.
We write $D_{s,c} :=D_s$ just for the remainder of this proof to stress the 
dependence on $c$.
The mapping $(c,u)\mapsto -D_{s,c}^2 u$ is smooth $U\x H^k_* \to H^{k-2n}$ and is
a bibounded linear isomorphism $H^k_*\to H^{k-2n}_*$ for fixed $c$. This can be
seen as follows (compare with \cite[4.5]{MM3}): 
It is true if $c$ is parametrized by arclength (look at it
in the space of Fourier coefficients). The index is invariant under
continuous deformations of elliptic operators of fixed degree, 
so the index of $-D_s^2$ is zero in general. But
$-D_s^2$ is self-adjoint positive, so it is injective with vanishing
index, thus surjective. By the open mapping theorem it is then bibounded.
Moreover $(c,w)\mapsto (-D_s^2)\i(w)$ is smooth $U^k\x H^{k-2n}_*\to H^k_*$ (by the
inverse function theorem on Banach spaces).
The mapping $(c,f)\mapsto D_s f = \frac1{|c_\th|}\p_\th f$
is smooth $H^k_*\x H^m_*\supset U\x H^m_*\to H^{m-1}$ for $k\ge m$, and is linear in $f$.
We have $v=D_{s,c} c$ and $n=i D_{s,c} c$. The mapping $c\mapsto \ka(c)$ is smooth
on the $H^2$-open set $\{c:|c_\th|>0\}\subset H^k_*$ into $H^{k-2}_*$.
Keeping all this in mind we now write the geodesic equation
\thetag{\ref{eq:geod}} as follows:

\begin{align*}
c_t &= u =:X_1(c,u)
\\
u_t &= -D_{s,c}^{-2}\Big(\tfrac12\|u|_{t=0}\|_G^2.\ka_c.n_c
  +\tfrac12D_{s,c}(|D_{s,c}u|^2.v_c)
  +D_{s,c}(\langle D_{s,c}u,v_c \rangle.D_{s,c}u)\Big)
\\&\quad
  -\Big( \frac1{\ell_c}\int\langle u,D_{s,c}^2 c \rangle\,ds\Big)\Big|_{t=0}\cdot u
\\&
=: X_2(c,u)
\end{align*}
Here we used that along any geodesic the norm $\|c_t\|_G$ and
the scaling momentum 
$- \frac1{\ell_c}\int\langle c_t,D_{s,c}^2 c \rangle\,ds=\p_t\log\ell(c)$
are both constant in $t$. 
Now a term by term investigation shows that 
the expression in the brackets is smooth $U^k\x H^k\to H^{k-2}$ since
$k-2>\frac12$. The operator $-D_{s,c}^{-2}$ then takes it smoothly back
to $H^k$. So the vector field $X=(X_1,X_2)$ is smooth on $U^k\x H^k$.
Thus the flow $\on{Fl}^k$ exists on $H^k$ and is smooth in $t$ and the initial 
conditions for fixed $k$.

Now we consider smooth initial conditions $c_0=c(0,\;.\;)$ and
$u_0=c_t(0,\;.\;)=u(0,\;.\;)$ in $C^\infty(S^1,\mathbb R^2)$. Suppose the 
trajectory $\on{Fl}^k_t(c_0,u_0)$ of
$X$ through these intial conditions in $H^k$ maximally exists for $t\in
(-a_k,b_k)$, and the trajectory $\on{Fl}^{k+1}_t(c_0,u_0)$ in $H^{k+1}$ 
maximally exists for $t\in(-a_{k+1},b_{k+1})$ with $b_{k+1}<b_k$. By
uniqueness we have $\on{Fl}^{k+1}_t(c_0,u_0)=\on{Fl}^{k}_t(c_0,u_0)$ for
$t\in (-a_{k+1,}b_{k+1})$. We now apply $\p_\th$ to the equation
$u_t=X_2(c,u)=-D_{s,c}^{-2}(\,\dots\,)$, 
note that the commutator $[\p_\th,-D_{s,c}^{-2}]$ is a pseudo differential
operator of order $-2$ again, and write $w=\p_\th u$. We obtain 
$w_t=\p_\th u_t= -D_{s,c}^{-2}\p_\th(\,\dots\,) +
[\p_\th,-D_{s,c}^{-2}](\,\dots\,) + \text{const}.w$.
In the term $\p_\th(\,\dots\,)$ we consider now only the terms $\p_\th^{3}u$ and
rename them $\p_\th^{2}w$.
Then we get an 
equation $w_t(t,\th)=\tilde X_2(t,w(t,\th))$ which is inhomogeneous bounded linear
in $w\in H^k$ with
coefficients bounded linear operators on $H^k$ 
which are $C^\infty$ functions of $c, u \in H^k$. 
These we already know on the intervall $(-a_k,b_k)$. 
This equation 
therefore has a solution $w(t,\;.\;)$ for all $t$ for which the
coefficients exists, thus for all $t\in (a_k,b_k)$. The limit 
$\lim_{t\nearrow b_{k+1}} w(t,\;.\;)$ exists in $H^k$ 
and by continuity it equals $\p_\th u$ in $H^k$
at $t=b_{k+1}$. Thus the $H^{k+1}$-flow was not
maximal and can be continued. So $(-a_{k+1},b_{k+1})=(a_k,b_k)$.  
We can iterate this and conclude that the flow of $X$ exists in
$\bigcap_{m\ge k} H^{m}= C^\infty$. 
\qed\end{demo}


\begin{thebibliography}{10}

\bibitem{alexhopf}
{\sc P.~Alexandroff, H~Hopf.} {\em Topologie. I}. 
Springer-Verlag, Berlin-New York, 1974

\bibitem{bcgj95}
{\sc R.~Basri, L.~Costa, D.~Geiger, and D.~Jacobs}, {\em Determining the
  similarity of deformable shapes}, in IEEE Workshop on Physics based Modeling
  in Computer Vision, 1995, pp.~135--143.

\bibitem{bl04} 
{\sc R.D. Benguriaand, M. Loss}, {\em Connection
between the Lieb-Thirring conjecture for Schr\"{o}\-din\-ger operators and
an isoperimetric problem for ovals in the plane'}, Contemporary
Math. {362} (2004), pp.53-61.

\bibitem{doc92} 
{\sc M.~P.~ Do Carmo}, {\em Riemannian Geometry} Mathematical Theory
and Applications, Birk\-h\"au\-ser, 1992.

\bibitem{gk91}
{\sc U.~Grenander and D.~M. Keenan}, {\em On the shape of plane images}, Siam
  J. Appl. Math., 53 (1991), pp.~1072--1094.

\bibitem{hel78}
{\sc S.~Helgason}, {\em Differential Geometry, Lie Groups and
Symmetric Spaces}, Graduate Studies in Mathematics, Vol. 34, American
Mathematical Society 1978 (revised version 2001).

\bibitem{ksm02}
{\sc E.~Klassen, A.~Srivastava, W.~Mio, and S.~Joshi}, 
{\em Analysis of planar   shapes using geodesic paths on shape spaces}, 
IEEE Trans. PAMI,  (2002).

\bibitem{klm04}
{\sc A.~Kriegl, M.~Losik, P.W.~Michor},
{\em Choosing roots of polynomials smoothly, II},  
Israel J. Math.  139 (2004), 183-188.

\bibitem{km93}
{\sc A.~Kriegl, P.W.~Michor}, 
{\em The Convenient Setting of Global Analysis},  
Mathematical Surveys and Monographs, Vol. 53, American Mathematical Society 1997.

\bibitem{ky06}
{\sc H.~Krim and A.~Yezzi}, eds., {\em Statistics and Analysis of Shapes},
Birkh\"auser, 2006.

\bibitem{MM3}
{\sc P.~Michor and D.~Mumford}, {\em An overview of the riemannian metrics on
  spaces of curves using the hamiltonian approach}, 
  Applied and Computational Harmonic Analysis, doi:10.1016/j.acha.2006.07.004.

\bibitem{msj05}
{\sc W.~Mio, A.~Srivastava, and S.~Joshi}, {\em On the shape of plane elastic
  curves}, tech. report, Florida State University, 1995.

\bibitem{ner01}
{\sc Y.~A. Neretin}, {\em On jordan angles and the triangle inequality in
  grassmann manifold}, Geometriae Dedicata, 86 (2001).

\bibitem{one66}
{\sc B.~O'Neill}, {\em The fundamental equations of a submersion}, Michigan
  Math. J.,  (1966), pp.~459--469.

\bibitem{sm04}
{\sc E.~Sharon and D.~Mumford}, {\em 2d-shape analysis using conformal
  mapping}, in Proceedings IEEE Conference on Computer Vision and Pattern
  Recognition, 2004.

\bibitem{ty00}
{\sc A.~Trouv\'e and L.~Younes}, {\em Diffeomorphic matching in 1d: designing
  and minimizing matching functionals}, in Proceedings of ECCV 2000, D.~Vernon,
  ed., 2000.

\bibitem{ty01}
{\sc A.~Trouv\'e and L.~Younes},   
 {\em On a class of optimal matching problems in 1 dimension},
 {Siam J. Control Opt.}, vol. 39, number 4, pp 1112-1135, 2001. 


\bibitem{you98}
{\sc L.~Younes}, {\em Computable elastic distances between shapes}, SIAM J.
  Appl. Math, 58 (1998), pp.~565--586.

\bibitem{you99}
{\sc L.~Younes}, {\em Optimal matching between shapes via elastic
  deformations}, Image and Vision Computing,  (1999).

\end{thebibliography}

\end{document}